\documentclass[11pt,final]{iopart}

\usepackage{amsfonts}
\usepackage{amsthm}
\usepackage{iopams} 
\usepackage{graphicx,subfigure}
\usepackage{caption}
\usepackage[dvips]{epsfig}
\usepackage{float}
\usepackage{showkeys}
\usepackage{todo}
\usepackage{algorithmicx}
\usepackage{algpseudocode}
\usepackage{prettyref}
\usepackage{verbatim} 
\usepackage[latin1]{inputenc}
\usepackage{tikz}
\usetikzlibrary{shapes,arrows}
\newrefformat{eq}{(\ref{#1})}
\newrefformat{fig}{Figure~\ref{#1}}
\usepackage[pdfpagelabels]{hyperref} 
\hypersetup{colorlinks=true,
linkcolor=blue,
citecolor=blue,
plainpages=false}
\newtheorem{theorem}{\sc Theorem}[section]

\newtheorem{definition}[theorem]{\sc Definition}

\newcommand{\mb}{\mathbf}
\newcommand{\argmin}{{\rm arg\,min}}

\eqnobysec

\makeindex

\newcommand{\cB}{{\mathcal B}}
\newcommand{\cC}{{\mathcal C}}
\newcommand{\cD}{{\mathcal D}}

\newcommand{\cG}{{\mathcal G}}
\newcommand{\cH}{{\mathcal H}}

\newcommand{\cL}{{\mathcal L}}

\newcommand{\cN}{{\mathcal N}}

\newcommand{\cR}{{\mathcal R}}

\newcommand{\cT}{{\mathcal T}}

\newcommand{\cV}{{\mathcal V}}
\newcommand{\cW}{{\mathcal W}}

\newcommand{\cZ}{{\mathcal Z}}

\newcommand{\bx}{{\bf x}}

\newcommand{\N}{\mathbb{N}}
\newcommand{\R}{\mathbb{R}}

\newcommand{\norm}[1]{|| #1||}

\renewcommand{\div}{{\rm div \;}}

\newcommand{\myarr}[1]{\left( \begin{array}{cccccccccccccccccc}#1\end{array}\right)}

\newcommand{\ba}{\begin{array}}
\newcommand{\ea}{\end{array}}
\newcommand{\be}{\begin{equation}}
\newcommand{\ee}{\end{equation}}
\newcommand{\bea}{\begin{eqnarray}}
\newcommand{\eea}{\end{eqnarray}}
\newcommand{\beq}{\begin{equation}}
\newcommand{\eeq}{\end{equation}}
\newcommand{\bqt}{\begin{quote}}
\newcommand{\eqt}{\end{quote}}

\begin{document}

%
\title[Refractivity Profile in Atmospheric Tomography]
{Quasi-Newton Approach for an Atmospheric Tomography Problem}

%
\author{ Erdem Altuntac }

\address{Institute for Numerical and Applied Mathematics,
University of G\"{o}ttingen, Lotzestr. 16-18,
D-37083, G\"{o}ttingen, Germany}


\ead{\mailto{e.altuntac@math.uni-goettingen.de}}

\begin{abstract}

This work studies the usage of well-known smoothed total variation regularization
for solving an atmospheric tomography problem named as {\em GPS-tomography}
in some quasi-Newton methods. That is we solve an unconstrained, convex, smooth
minimization problem associated with a general type Tikhonov functional 
containing smoothed form of total variation penalty term
by quasi-Newton methods. As a result of the conducted experiments, 
on the basis of error analysis {\em i.e.} convergence analysis,
it is concluded that the limited memory BFGS algorithm with trust region 
is the most effective algorithm in terms obtaining a reasonable optimum solution.

\end{abstract}

\bigskip


\section{Introduction}
\label{intro}

One important predictor in meteorology is the humidity of the
atmosphere. This is estimated by fan-beam measurements between 
satellite transmitters and land-based receivers.  The measurements 
are sparse and fluctuate randomly with receiver availability. 
The task is to reconstruct from these measurements the 3-dimensional, 
spatially varying index of refraction of the atmosphere, 
from which the relative humidity can be inferred.

GPS-tomography involves the reconstruction of some quantity, pointwise within a volume 
(e.g. humidity) from geodesic X-ray measurements transmitted by nonuniformly distributed
transducers (satellites). These measurements are collected by nonuniformly
distributed receivers on the ground (ground stations).
As with conventional tomography, the task here is the reconstruction of the density 
volume profile of a layer in the atmosphere from a set of line integrals.
Function reconstruction from its measured line integrals 
was firstly proposed and solved in \textbf{\cite{Radon17}}. 
Profound mathematical and numerical aspects of the computerized tomography 
have been studied in \textbf{\cite{Natterer01, NattererWuebbeling01}}. 
Measurement from the Radon transform is obtained by integrating
some integrable function over the hyperplanes in $\R^{N}.$ The ray transform,
on the other hand, produces measurement by integrating the function over straight lines.
It is known that in the two dimensional tomography, general Radon and ray transformations 
coincide, \textbf{\cite[p. 17]{NattererWuebbeling01}}.

In the discretized form of the problem, it is assumed that each station receives
equal number of signals transmitted by the satellites. Also for the sake of 
simplicity, we ignore any deviations from the shortest path between transmitters
and receivers due to atmospheric refractivity. The received signal is then
modelled as a line integral along the shortest path between the satellites
and the ground stations.

Peculiar to this problem, reconstructions by Kalman filtering
and ART have been widely applied, 
\textbf{\cite{Bender11, MiidlaRannatUba08, Perler11, ZusBender12}}.
Different from these conventional numerical reconstruction methods, we 
propose a quasi-Newton approach. One of the effective quasi-Newton 
methods is {\em limited memory BFGS} (L-BFGS) algorithm which
is particularly suggested by this work.
The L-BFGS algorithm has been also applied for atmospheric imaging 
wherby the forward problem has been modelled as a phase retrieval problem, see \textbf{\cite{Vogel00}}.
We, on the other hand, consider the forward model as a linear atmospheric transmission problem 
which is a {\em straight line approximation}.
This means that despite the refractivity in the microwave signals while traversing the troposphere
layer of the atmosphere, we ignore attenuation. 
The unknown function is denoted by $\varphi$ which is assumed
to be in the class of some reflexive Banach space $\cV = \cL^{p}(\Omega),$
for $1 \leq p \leq d/(d-1)$ where $d = 3$ since this work 
focuses on three dimensional reconstruction. The measured noisy
data is assumed to be in the class of some Hilbert space $\cH.$ 
We, then, seek the minimizer for some general Tikhonov objective 
functional given in the form of

\bea
\label{cost_functional}
F_{\alpha}(\varphi , f^{\delta}) : & \cV \times \cH & \longrightarrow \R_{+}
\nonumber\\
& (\varphi , f^{\delta}) & \longmapsto F_{\alpha}(\varphi , f^{\delta}) := 
\frac{1}{2}\norm{\cT\varphi - f^{\delta}}_{\cH}^2 + \alpha J(\varphi) ,
\eea
where the forward operator $\cT :\cV \rightarrow \cH$, as will
be described soon, is a linear fan-beam projection operator. Here, 
the penalty term $J :\cV \rightarrow \R_{+}$ is convex and 
Fr\'{e}chet differentiable with the regularization parameter 
$\alpha > 0$ before it.

We demonstrate our regularization on simulated data, employing a novel
reverse-communication large-scale nonlinear optimization software 
SAMSARA which has been developed by D. R. Luke \textbf{\cite{Luke}}.
Comparison between the illustrated results from SAMSARA and 
the results from traditional lagged diffusivity fixed point
iteration algorithm, {\em LDFP} in \textbf{\cite{Vogel02, VogelOman96}}, 
is also provided.


\subsection{Physical Problem: From Propagation in Time to Propagation in Space}
\label{physical_problem}


This is an inverse problem with incomplete data.
It is well known that the incompleteness of data causes 
nonuniqueness issue in inverse problems,
\textbf{\cite[p. 144]{NattererWuebbeling01}}.
Particularly in tomography, the assumption of compact support is essential 
in order for unique solvabilibility. In other words,
problems characterizing incomplete data case are uniquely 
solvable if the unknown function $\varphi$ has compact support.
In this subsection, although it does not completely overlap
the reality, we will model the physical problem
with the geometrical assumption of compact support.
Firstly, just by the nature of the physical problem, $\varphi$ 
is not a constant function and contains smooth intensity. 
Formulation of the simulated profile
is presented in Subsection \ref{the_synthetic_profile}.
Since this work solely aims to provide empirical results 
for a large scale application problem by some well-known
optimizatin and regularization strategies, we will not state
any theoretical result. However, still as a duty of any inverse problem research work,
formal assumptional statemens on compact support must be made that uniqueness
principle is verifiable.

\begin{figure}[ht]
\centering
\begin{minipage}{.7\textwidth}
  \centering

\includegraphics[width=0.65\textwidth]{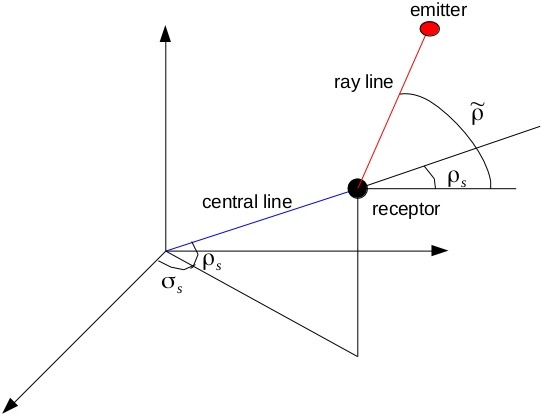}
\caption[Sketch of the angular parameterization of the tomography problem.]
{\footnotesize The sketch of angular parameterization. Intersection point
between the central and the ray lines is the receptor with
the angles $(\rho_s , \sigma_s).$ Emitter is located at the highest altitude $h^{\infty}.$
Its elevation angle according to the surface is denoted by $\tilde{\rho}.$}
\label{angular_parameterization}
\end{minipage}\hfill
\begin{minipage}{.7\textwidth}
  \centering
\includegraphics[width=1\textwidth]{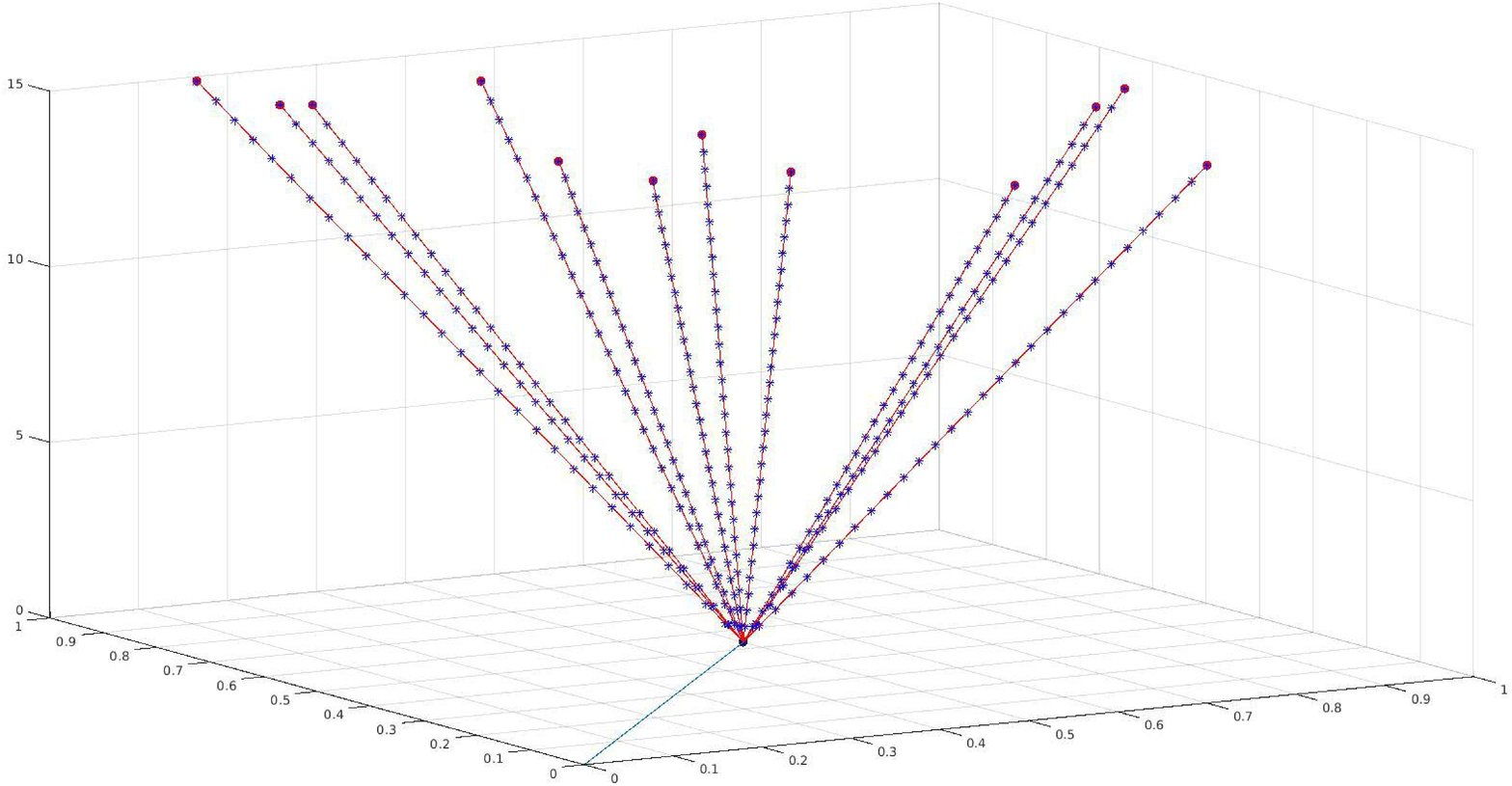}
\caption{\footnotesize Finitely number of points denoted by blue stars along $12$ rays are illustrated.
These points are found via the parameter function $t(\epsilon) = \frac{\epsilon}{\sin{(\tilde{\rho})}},$
for $\epsilon \in [z_s , h^{\infty}]$ where $h^{\infty}$
is the upper bound of the line integral in (\ref{integral_transform1}).}
\label{intro_3D_Network_int_points}
\end{minipage}
\end{figure}

Let $(\rho_s , \sigma_s)$ be the polar angles 
of the station $\mathbf{s}$ as inclination and azimuth respectively.
Then in spherical coordinates, the location of any station $\mathbf{s}$ 
is given by
\begin{displaymath}
\mathbf{s} = R \myarr{ \cos(\rho_s)\cos(\sigma_s), \cos(\rho_s)\sin(\sigma_s), \sin(\rho_s) }^{T},
\end{displaymath}
where $ R = \sqrt{{x_s}^2 + {y_s}^2 + {z_s}^2}.$
Following \textbf{\cite[Ch. 2]{Natterer01}} and \textbf{\cite[p. 45]{OlafssonQuinto06}}, 
the signal path direction $\vec{\theta}$ is reparametrized by 
\bea
\vec{\theta} & = &
-\myarr{ \cos(\rho_r + \rho_s)\cos(\sigma_r + \sigma_s), \cos(\rho_r + \rho_s)\sin(\sigma_r + \sigma_s), \sin(\rho_r + \rho_s)}^{T} 
\nonumber\\
& = & - \myarr{ \cos(\tilde{\rho})\cos(\tilde{\sigma}), \cos(\tilde{\rho})\sin(\tilde{\sigma}) , \sin(\tilde{\rho}) }^{T} ,
\eea
where the inclination and the azimuth of the signal path according to
the surface are denoted by 
$(\rho_r + \rho_s , \sigma_r + \sigma_s) = (\tilde{\rho} , \tilde{\sigma}),$
see Figure \ref{angular_parameterization} for this angular parameterization. 


Let $g$ be some Lipschitz continuous 
function with its Lipschitz constant $L_g \in \R_{+}$ 
for the surface of the earth,
\bea
\label{surface_function}
g : [0, S] \times [0, P] \rightarrow [0 , h^{\infty}) ,
\eea
and denote by $\cG$ the graph of the surface function $g$
\bea
\label{surface_graph}
\cG = \mathrm{graph}(g)
:= \{(x, y, z) \mbox{ } \vert \mbox{ } (x , y) \in [0, S] \times [0, P] , z = g(x, y) \geq 0 \} .
\eea
A ground station $\mb{s}$ is a set of points in $\R^{3}$
located on earth with the coordinate points $(x_s, y_s, z_s),$
\bea
\label{ground_station}
\mb{s} := (x_s, y_s, z_s) \in \cG , 
\eea
and likewise emitters $\mb{e}$ that are all located 
at the same altitude $h^{\infty}$ is also set of 
points in $\R^{3},$
\bea
\label{emitter}
\mb{e} := \{ (x_e, y_e, z_e) \mbox{ } \vert \mbox{ } (x_e, y_e) \in [0, S] \times [0, P] \mbox{ and } z_e = h^{\infty} \} .
\eea
Our area of interest is a compact subdomain,
{\em i.e.} $\Omega \subset \Omega_{o},$
\bea
\label{omega}
\Omega_{o} := \left\{ (x, y, z) \in  [0, S] \times [0, P] \times [0, h^{\infty}]\mbox{ } \vert \mbox{ } z \geq g(x, y) \right\} .
\eea
Obviously, the definitions in (\ref{surface_graph}) and (\ref{omega}) both imply
$\cG \subset \Omega_{o}.$


Since we consider our network as straight line approximation,
that is we do not include attenuation,
we model each signal path as a ray in $\R^{3}.$
There can be formulated a linear parameter function 
$t : \R \rightarrow \R,$ $t(\epsilon) := \frac{\epsilon}{\sin{(\tilde{\rho})}},$ 
such that a ray in $\R^{3}$ starting from the station $\mb{s}$ 
in the direction $\vec{\theta} \in \mathbb{S}^2$ 
is defined by 
\beq
\label{signal_set}
\gamma_{[\mathbf{s} , \vec{\theta}]}(\epsilon) := \left\{ \mathbf{s} + \frac{\epsilon}{\sin{(\tilde{\rho})}} \vec{\theta} \mbox{ } 
\Bigg\vert \mbox{ } \epsilon \in [z_s , h^{\infty}] \right\}.
\eeq
Here, in fact, $\gamma$ is the minimal path between any two points
in $\R^{3}.$ So, a microwave signal takes the least time 
$T$ with speed $c$ along this path
\bea
\label{fermats_principle}
T = \frac{1}{c} \int_{\bf{r} \in \gamma} {n(\bf{r})} d \bf{r} ,
\eea
where $n$ is index of refraction.
The linear relation between the refractivity profile $N$
and the refractive index $n$ is expressed by $N = 10^{6}(n-1),$
\textbf{\cite{Bender11, MiidlaRannatUba08, ZusBender12}}.
Thus, if one chooses the refractivity profile as the frame of reference,
then (\ref{fermats_principle}) reads
\bea
\label{fermats_principle2}
T = \frac{1}{c} \int_{\bf{r} \in \gamma} (10^{-6}{N(\bf{r})} + 1) d \bf{r} .
\eea
To obtain measurement $f,$ we apply fan-beam projection operator along the
ray $\gamma_{[\mathbf{s} , \vec{\theta}]}$ on some density profile defined by 
$\varphi := \frac{1}{c}\left( 10^{-6}N + 1 \right) = \frac{n}{c}.$
The unknown density function $\varphi : \Omega \rightarrow \R$
is assumed to be integrable and, by convention, vanishes outside the area of interest $\Omega.$
This is explained by introducing a step function as such

\bea
\tilde{\varphi}(x) :=
\left\{ \begin{array}{rcl}
\varphi(x) & \mbox{, for}
& x \in \Omega \\ 
0 & \mbox{, for} & x \in \Omega_{o}\backslash\Omega .
\end{array}\right . 
\eea
Physically, there exist many rays in various directions $\vec{\theta} \in \mathbb{S}^2.$
However, the measured data can only be obtained through the rays which do not have empty 
intersection with the area of interest $\Omega .$ Let $\cZ$ be the domain of the integrated
measurement which is the function of station $\mathbf{s}$ and directional vector $\vec{\theta}.$ 
Denote by
\bea
\label{direction_set}
\mathbb{S}_{\mathbf{s}} := \{ \vec{\theta} \in \mathbb{S}^2  \vert (\mathbf{s} , \vec{\theta}) \in \cZ \} ,
\eea
the set of intercepted directions where the domain of the integrated measurement 
through one ray $\gamma_{[\mb{s},\vec{\theta}]}$ can be presented by
\bea
\label{measurement}
f : \cZ = \cD(f) \subset \cG \times \mathbb{S}_{\mathbf{s}} \rightarrow \R_{+},
\eea
with
\bea
\label{measurement_domain}
\cD(f) := \{(\mb{s} , \vec{\theta}) \mbox{ } \vert \mbox{ } \tilde{\rho} \geq \vert \arctan(L_g) \vert , \pi - \tilde{\rho} > 0 
\mbox{ , and }  \gamma_{[\mathbf{s} , \vec{\theta}]} \cap \Omega \neq \emptyset \} .
\eea
By (\ref{measurement_domain}), one must understand that the slope
of the ray cannot be larger than the elevation angle $\tilde{\rho}.$
Furthermore, rays that are parallel to the surface are not taken into account
for the measurement. There could also be rays 
that do not intersect with the area of interest $\Omega.$
Therefore, we are only interested in the rays that have
no empty intersection with $\Omega,$

\begin{displaymath}
\gamma_{[\mb{s} , \vec{\theta}]} \cap \Omega \neq \emptyset \mbox{, for } \theta \in \mathbb{S}_{\mb{s}} \subset \mathbb{S}^2 
\mbox{, }  \mb{s} \in \cG \subset \Omega_{o}.
\end{displaymath}
Then, in fact, the measured data 
$f({\mathbf{s}},\vec{\theta})$
is obtained only for $\vec{\theta} \in \mathbb{S}_{\mathbf{s}} .$
Note that $\mathbb{S}_{\mathbf{s}} \subset \mathbb{S}^2,$ which is
the partial information case. Thus collection of the measurement 
operation, in light of fan-beam projection principle, is formulated by
\beq
\label{integral_transform}
f(\mb{s} , \vec{\theta}) = \cT_{\mb{s}}\varphi(\vec{\theta}) = (\cT\varphi)(\mathbf{s} , \vec{\theta}) = \int_{\gamma_{[\mathbf{s} , \vec{\theta}]}} \varphi(\mb{r}) d \mb{r} , 
\mbox{ for } \mb{r} \in \gamma_{[\mathbf{s} , \vec{\theta}]}\subset \Omega_{o} .
\eeq
Also, with the angular parameterization introduced above,
we then have
\begin{displaymath}
\cT_{\mb{s}}\varphi(\vec{\theta}) = \cT\varphi(\tilde{\rho} , \tilde{\sigma}), \mbox{ for }
(\tilde{\rho} , \tilde{\sigma}) \in (0 , \pi) \times (0 , 2\pi) .
\end{displaymath}
According to \textbf{\cite[Theorems 5.1 - 5.6]{HamakerSmith80}} and 
\textbf{\cite[Theorem 6.2]{NattererWuebbeling01}}
the linear transformation (\ref{integral_transform}) is injective 
only under compact support assumption and in the presence of
directional vectors from the set of intercepted directions,
$\theta \in \mathbb{S}_{\mathbf{s}}.$ It is not possible
to reconstruct the unknown function $\varphi$ exactly from
finitely number of measurements. However, 
\textbf{\cite[Theorems 5.1 - 5.6]{HamakerSmith80}}
show that arbitrarily good approximation can be obtained.


The discretized integration from one point
to the next one along the ray $\gamma$
is carried out via the parameter function 
$t(\epsilon) = \frac{\epsilon}{\sin{(\tilde{\rho})}},$
for any $\epsilon \in [z_s , h^{\infty}],$ 
see Figure \ref{intro_3D_Network_int_points}.
In the continuum form,
we use ray transform in the direction
$\vec{\theta}(\tilde{\rho} , \tilde{\sigma}) \in \mathbb{S}_{\mathbf{s}}$
for any angle pairs $(\tilde{\rho} , \tilde{\sigma}),$
on the density function $\varphi : \Omega \subset \R^{3} \rightarrow \R_{+}$ 
as such

\bea
\label{integral_transform1}
f(\mathbf{s}, \vec{\theta}) = \cT_{\mathbf{s}}\varphi(\vec{\theta}) = (\cT \varphi)(\mathbf{s} , \vec{\theta})
= \int_{\mathbf{r} \in \gamma_{[\mathbf{s}, \vec{\theta}]}} \varphi({\mathbf{r}}) d \mathbf{r}
= \int_{z_s}^{h^{\infty}} \varphi(\gamma_{[\mathbf{s}, \vec{\theta}]}(\epsilon)) \vert \gamma^{\prime}_{[\mathbf{s}, \vec{\theta}]}(\epsilon) \vert d \epsilon ,
\eea
where 
\bea
\label{linearity_line_int}
\vert \gamma^{\prime}_{[\mathbf{s}, \vec{\theta}]}(\epsilon) \vert
= \left\vert \frac{\vec{\theta}}{\sin{(\tilde{\rho})}} \right\vert, \mbox{ with } \vert \vec{\theta} \vert = 1.
\eea
The representation (\ref{integral_transform1}) is comparable with 
its nonlinear counterpart in \textbf{\cite[Eq. (1.3)]{SchroederSchuster15}}. 
So as a linear operator equation, we have $\cT \varphi = f$
where $\cT$ represents the line integration operating on
the density profile $\varphi$ to obtain measurement $f.$




\begin{figure}[ht]
  \centering
    \includegraphics[height=4.5in,width=6.3in,angle=0]{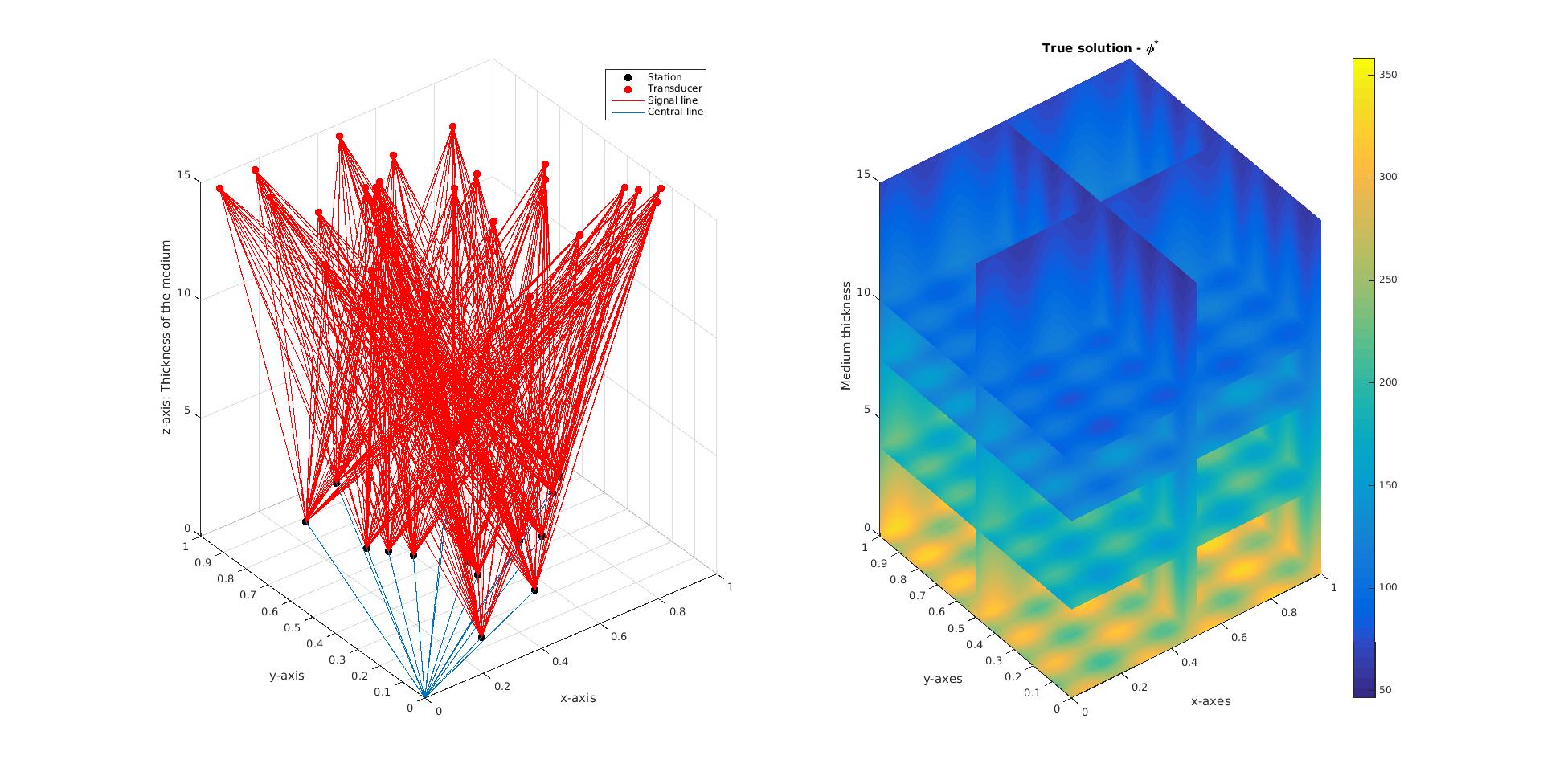}
\caption[A 3-D network together with simulated data is illustrated over a non-uniformly scaled domain.]
{\footnotesize A 3-D network together with simulated data $\varphi,$
as the true solution $\varphi^{\dagger},$
is illustrated over a non-uniformly scaled domain.
Black dots indicate stations whilst signals penetrate the area of interest through red dots.
In this illustration, $15$ ground stations (receiver) intercept signals emitted by $30$ 
transducers and all are randomly distributed over $27000$ grid nodes.}
\label{3-D_Network_numerics_intro}
\end{figure}

\section{Minimization Problem, Existence and Uniqueness of the Regularized Solution}
\label{minimization_problem}

It has been conveyed that the use of $TV$
promotes sparsity of the gradient, \textbf{\cite{BenningGladdenHolland14}}.
In our numerical illustrations,
we have simulated a data with smooth intensity, see Subsection \ref{the_synthetic_profile}.
The weak formulation of TV of some function $\varphi$ defined over the compact domain $\Omega$ is given below.

\begin{definition}\textbf{[$TV(\varphi , \Omega)$]}\textbf{\cite[Definition 9.64]{ScherzerGrasmair09}}
\label{TV_functional}
Over the compact domain $\Omega,$
total variation of a function $TV(\varphi , \Omega)$ 
is defined in the weak sense as follows ,

\bea
\label{TV_functional2}
 TV(\varphi , \Omega) := \sup_{\Phi\in \cC_{c}^{1}(\Omega)} 
\left\{ \int_{\Omega} \varphi(x) \div \Phi(x) dx\mbox{ } : \mbox{ } \norm{\Phi}_{\infty} \leq 1 \right\} .
\eea
\end{definition}
Total variation type regularization targets 
the reconstruction of bounded variation (BV) 
class of functions that are defined by
\bea
\label{bv_def}
BV(\Omega) := \{ \varphi \in \cL^{1}(\Omega) : TV(\varphi , \Omega) < \infty \},
\eea
endowed with the norm
\beq
\norm{\varphi}_{BV} := \norm{\varphi}_{\cL^1} + TV(\varphi , \Omega).
\eeq
BV function spaces are Banach spaces, \textbf{\cite{Vogel02}}.
By the result in \textbf{\cite[Theorem 2.1]{AcarVogel94}}, 
it is known that one can arrive, with a proper choice of $\Phi \in \cC_{c}^{1}(\Omega),$ 
in the following from (\ref{TV_functional2}),

\beq
\label{tv_integral_form}
TV(\varphi) = \int_{\Omega} \vert\nabla\varphi(x)\vert dx \cong
\int_{\Omega} \left( \vert\nabla\varphi(x)\vert^2 + \beta \right)^{1/2} dx ,
\eeq
where $0 < \beta < 1$ is fixed and the classical Euclidean norm 
is denoted by $\vert \cdot \vert.$ We also refer
\textbf{\cite{ChambolleLions97, ChanGolubMulet99, DobsonScherzer96, RudinOsherFatemi92, VogelOman96}}
where the smoothed form of (\ref{tv_integral_form}) has appeared.

With this theoretical motivation having stated,
we are tasked with constructing the regularized solution
$\varphi_{\alpha}^{\delta}$ over some compact 
and convex domain $\Omega \subset \cR^{3}$ by solving the following 
smooth, unconstrained, minimization problem,

\beq
\label{problem}
\varphi_{\alpha}^{\delta} \in \argmin_{\varphi \in \cV} 
\left\{ \frac{1}{2} \norm{\cT \varphi - f^{\delta}}_{\cH}^2 + \alpha J_{\beta}^{\mathrm{TV}}(\varphi) \right\} ,
\eeq
with its regularization parameter $\alpha > 0$ and
for the penalty term 
$J_{\beta}^{\mathrm{TV}} : \cV \rightarrow \R_{+},$ 
where in particular $\cV = \cW^{1,p}$ for $1 \leq p \leq d/(d-1)$
and $d =3,$ defined by

\bea
\label{smoothed_tv_regul}
J_{\beta}^{\mathrm{TV}}(\varphi) := \int_{\Omega} \sqrt{\vert\nabla\varphi(x)\vert^2 + \beta} d \bx .
\eea
It is the obvious property of the chosen penalty term
that $J_{\beta}^{\mathrm{TV}} \in C^{\infty}(\cV).$
Existence and uniqueness of the solution $\varphi_{\alpha}^{\delta}$ 
for the problem (\ref{problem}) has been studied extensively in 
\textbf{\cite{AcarVogel94}}. By the given facts of our forward operator,
one of which is that there could be rays with empty intersection,
it can be stated that

\bea
\label{property_T}
\cT 1 \neq 0 \mbox{, and } \cT 1 \geq 0.
\eea
This implies the $BV$ coercivity of the objective functional $F_{\alpha}$
from which the existence of the regularized solution is guaranteed.
Uniqueness of the solution is simply the consequence of the 
strict convexity $F_{\alpha}$ which is implied by the injectivity
of the forward operator $\cT.$


\section{Discretized Form of the Minimization Problem and the Toy Model Setup}
\label{discretized_minimization}

In the computerized environment we always work 
with finite dimensional setup, thus we only collect discrete data. 
So, we now introduce our tomographic application and the minimization 
problems with their components in the finite dimension. We consider the domain
$\Omega = [0, 1] \times [0, 1] \times [0, 15]$
and the meshsize $\Delta_{\bx} = 1/(N - 1)$
with some determined mesh point number $N \in \N$
for any point $\bx = (x, y, z) \in \Omega \subset \R^{3}.$
Note that, here $h^{\infty} = 15$ according to
(\ref{integral_transform1}).
Within our compact domain $\Omega \subset \R^{3},$
we then generate a point-to-point discretization 
by starting from some point $\bx_{i-1} \in \Omega \subset \R^{3}$ 
and iterating onward as such

\begin{displaymath}
\bx_{i} = \bx_{i - 1} + \Delta_{\bx}, \mbox{ for each } i = 2, \cdots, N .
\end{displaymath}
%
In our experiments, we have developed 
$N_x = N_y = N_z = 30$ nodes.
In the toy model, the speed of light is taken as $c = 1,$ 
see (\ref{fermats_principle}), in order to be able to
measure the propagation of the light beams in space instead of in time. 
Recall from the Section \ref{physical_problem}
by (\ref{integral_transform1}) that the electromagnetic 
signals with the angles $(\rho_r, \sigma_r)$ 
arrive in any receiver $\mathbf{s}$ with the polar angles 
$(\rho_s , \sigma_s)$ in various directions 
$\mathbf{\theta} \in \mathbb{S}^2.$ 
So the ray path in $\R^{3}$ is the set

\beq
\label{signal_set_discretized}
\gamma_{[\mathbf{s},\theta]}(\epsilon) := \left\{ \mathbf{s} + \frac{\epsilon}{\sin{(\tilde{\rho})}} \theta \hspace{.1cm} \Bigg\vert \hspace{.1cm} 
\epsilon \in [z_s , h^{\infty}] \right\} ,
\eeq
and the integral transformation that is used for data collection

\bea
f(\mathbf{s}, \theta) = \cT_{\mathbf{s}}\varphi(\theta) = (\cT \varphi)(\theta , \mathbf{s})
= \int_{\mathbf{r} \in \gamma_{[\mathbf{s}, \theta]}} \varphi({\bf{r}}) d \mathbf{r}
= \int_{z_s}^{h^{\infty}} \varphi(\gamma_{[\mathbf{s}, \vec{\theta}]}(\epsilon)) \vert \gamma^{\prime}_{[\mathbf{s}, \vec{\theta}]}(\epsilon) \vert d \epsilon .
\nonumber
\eea
The full path of the signal is the sum of the paths in the intercepted
grid nodes. The model can be interpreted as a system of linear equations. 
Let us denote the discretized integration by $T.$
With additive white Gaussian noise model vector 
$z_j \sim \cN(0, 1)$ (cf. \textbf{\cite{KekkonenLassasSiltanen14}}) 
and some known noise level $\delta,$ 
we produce measurement vector by

\bea
\label{disc_integration}
[T\varphi]_{j} = \sum_{i = 1}^{N} \varphi_{i} w_{i,j} = f_{j}^{\dagger} + \delta z_j = f_{j}^{\delta},
\eea
where $j = 1, 2, \cdots S,$ $S$ is the total number of signal
paths from all visible satellites in the network at a fixed
time instant, $N$ is the total number of grid nodes, $w_{i,j}$
is the length of $j^{th}$ ray passing through the node $i,$ 
the $\varphi_{i}$ is interpreted as the density of the corresponding
$i$th node, \textbf{\cite{MiidlaRannatUba08}}.

The parameter function $t(\epsilon) = \frac{\epsilon}{\sin{(\tilde{\rho})}}$
in (\ref{signal_set_discretized})
permits one to determine the points along each signal for 
any $\epsilon \in [z_s , h^{\infty}]$ where $h^{\infty}$
is the upper boundary of the medium as well as
the line integral in (\ref{integral_transform1}), 
see Figure \ref{intro_3D_Network_int_points}.

Regarding the discretized form of our minimization
problem (\ref{problem}) with its components, we are provided
with the compact forward operator 
$T : \R^{N} \rightarrow \R^{M}$ and 
the measurement vector $f^{\delta} \in \R^{M}.$
With this information, our cost functional is then 
$F_{\alpha}(\varphi , f^{\delta}) : \R^{M \times N} \rightarrow \R_{+},$ 
and we seek for the optimum solution to the problem

\bea
\label{discretized_problem}
\varphi_{\alpha}^{\nu + 1} \in \argmin_{\varphi^{\nu} \in \R^{N}} \left\{ 
\frac{1}{2}\norm{T\varphi^{\nu} - f^{\delta}}_{2}^{2} + \alpha J(\varphi^{\nu}) \right\} .
\eea
Since we have focused on the smoothed form of the total variation
regularization in our analysis, we then define the smooth-TV
penalty by
\beq
\label{discretized_TV}
J_{\beta}^{\mathrm{TV}}(\varphi^{\nu}) := \sum_{i=1}^{n_x}\sum_{j=1}^{n_y}\sum_{k=1}^{n_z}
\Gamma_{\beta} \left( (D_{ijk}^x \varphi^{\nu})^2 + (D_{ijk}^y \varphi^{\nu})^2 + (D_{ijk}^z \varphi^{\nu})^2 \right) \Delta_{x} \Delta_{y} \Delta_{z},
\eeq
where the smoothing functional 
$\Gamma_{\beta}(\Phi) := \sqrt{\vert \Phi \vert^2 + \beta}$
for some fixed $\beta \in (0,1)$ and the discretized spatial derivatives 
according to the central difference form

\beq
\label{discrete_derivative}
D_{ijk}^{x}\varphi = \frac{\varphi_{i+1,j,k}^{\nu} - \varphi_{i-1,j,k}^{\nu}}{2\Delta_{x}}, \hspace{.3cm}
D_{ijk}^{y}\varphi = \frac{\varphi_{i,j+1,k}^{\nu} - \varphi_{i,j-1,k}^{\nu}}{2\Delta_{y}}, \hspace{.3cm}
D_{ijk}^{z}\varphi = \frac{\varphi_{i,j,k+1}^{\nu} - \varphi_{i,j,k-1}^{\nu}}{2\Delta_{z}}.
\eeq
The optimum solution $\varphi_{\alpha}^{\nu}$ 
must satisfy the first optimality condition. That is

\begin{displaymath}
0 = \nabla F_{\alpha}(\varphi_{\alpha}^{\nu} , f^{\delta}) = 
T^{\ast}(T\varphi_{\alpha}^{\nu} - f^{\delta})
+ \alpha \nabla J_{\beta}^{\mathrm{TV}}(\varphi_{\alpha}^{\nu}) .
\end{displaymath}
Here $\nabla J_{\beta}^{\mathrm{TV}}(\varphi)$ is calculated by 
$\frac{d}{dt} J_{\beta}^{\mathrm{TV}}(\varphi + t \Psi) \vert_{t = 0}$ in the direction
$\Psi \in \cC_{c}^{1}(\Omega)$ such that $\norm{\Psi} \leq 1.$
It can be observed that 
$\nabla J_{\beta}^{\mathrm{TV}}(\varphi) = L(\varphi)\varphi$ 
with the nonlinear term $L(\varphi),$

\bea
\label{gradient_varphi}
L(\varphi) & = & D_{x}^T diag(\Gamma^{'}(\varphi)) D_{x} + D_{y}^T diag(\Gamma^{'}(\varphi)) D_{y} + D_{z}^T diag(\Gamma^{'}(\varphi)) D_{z}
\nonumber\\
& = & \myarr{D_{x}^T & D_{y}^T & D_{z}^T} 
\myarr{diag(\Gamma^{'}(\varphi)) & 0 & 0 \\ 0 & diag(\Gamma^{'}(\varphi)) & 0 \\0 & 0 & diag(\Gamma^{'}(\varphi))}
\myarr{D_{x} \\ D_{y} \\ D_{z}} .
\nonumber
\eea



\begin{figure}[ht]
\centering
\includegraphics[height=4.5in,width=5.5in,angle=0]{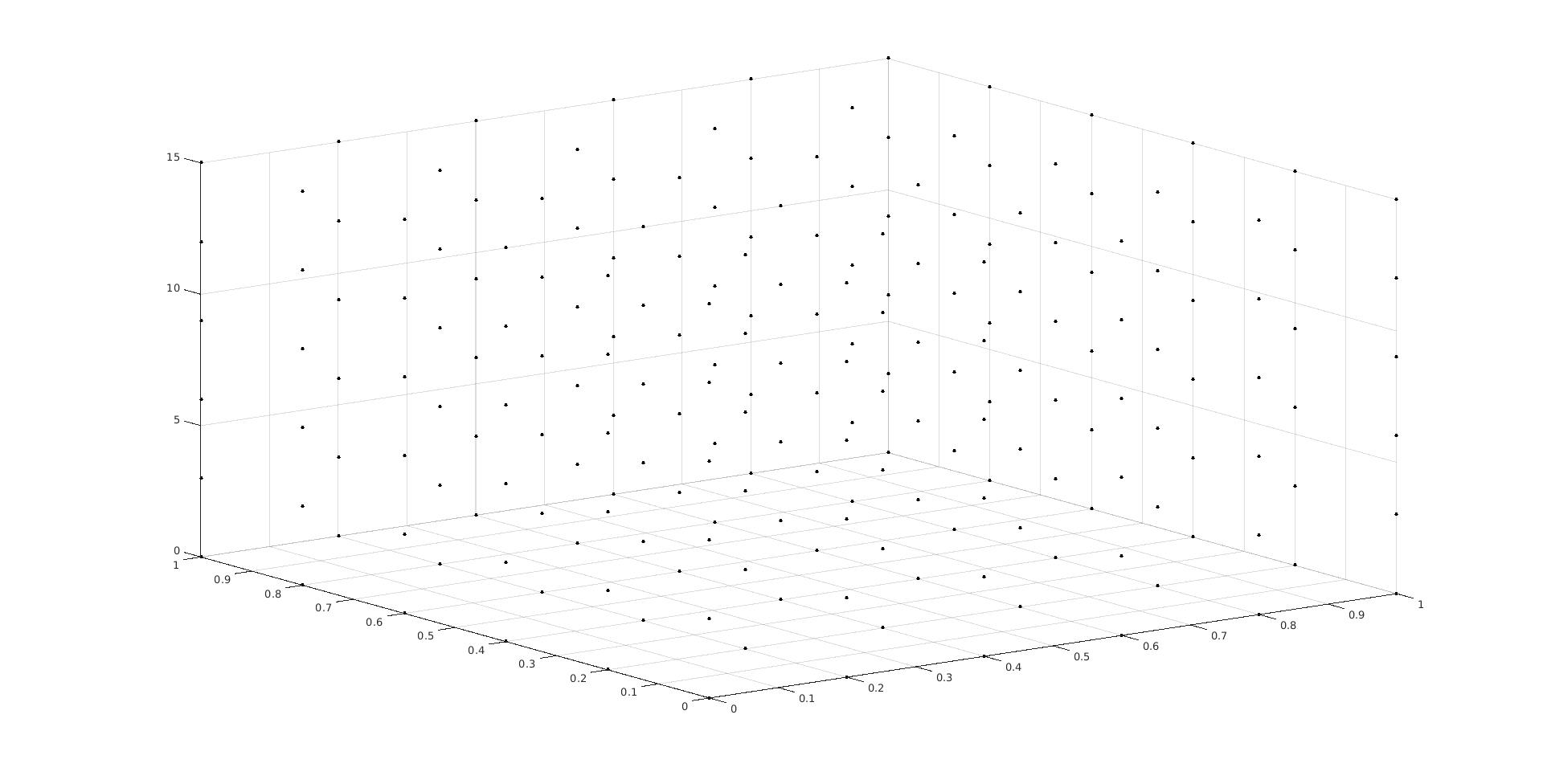}
\caption[Discretization of our area of interest.]
{\footnotesize Discretization of our area of interest. For illustration purpose, we only
present $6 \times 6 \times 6$ grid nodes.}
\label{3D_Network_pixel_points}
\end{figure}

\subsection{Empirical convergence analysis}
\label{empirical_convergence}

Recall that we aim to obtain approximate regularized solution
by solving the unconstrained, smooth minimization problem
\bea
\varphi_{\alpha}^{\delta} \in \argmin_{\varphi \in \cV}
\frac{1}{2} \norm{\cT \varphi - f^{\delta}}_{\cH}^2 + \alpha J_{\beta}^{\mathrm{TV}}(\varphi) ,
\nonumber
\eea
with the smooth-TV penalty term
\bea
J_{\beta}^{\mathrm{TV}}(\varphi) := \int_{\Omega} \sqrt{\vert\nabla\varphi(x)\vert^2 + \beta} d x .
\nonumber
\eea
Thus, we must observe sufficient decay
in the following components that we can claim the
optimum solution as a result of any algorithm;

\begin{itemize}
\item $F(\varphi_{\alpha}^{\nu} , f^{\delta});$ the functional value at every
updated point $\varphi_{\alpha}^{\nu},$

\item $\norm{\varphi_{\alpha}^{\nu} - \varphi_{\alpha}^{\nu - 1}};$ norm
of the succesive iterations at each iteration step $\nu = 1, 2, \cdots ,$

\item $\frac{\norm{\varphi_{\alpha}^{\nu} - \varphi^{\dagger}}}{\norm{\varphi^{\dagger}}};$
the relative error value of the reconstruction against the true solution $\varphi^{\dagger},$

\item $\Vert\nabla F(\varphi_{\alpha}^{\nu})\Vert;$ the norm of the gradient value of the functional
at every updated point $\varphi_{\alpha}^{\nu},$

\item $\Vert \cT\varphi_{\alpha}^{\nu} - f^{\delta} \Vert;$ the discrepancy of the image of the solution
against the given data $f^{\delta}.$

\end{itemize}

It is expected from the chosen regularization strategy that
this strategy must produce a reliable regularized minimizer 
$\varphi_{\alpha}^{\delta}.$ This reliability is tested in the 
framework of convergence concept. In order to be able to speak about 
the convergence of the regularized minimizer (the solution) 
$\varphi_{\alpha}^{\delta},$ there must be some reference solution 
to which the regularized solution will approximately converge during 
the iteration. Likewise in many inverse problems research works, 
we choose our reference solution as the true solution $\varphi^{\dagger}.$ 
Convergence of the regularized solution $\varphi_{\alpha}^{\delta}$ to the true
solution $\varphi^{\dagger}$ in the Hilbert norm sense 
$\norm{\varphi_{\alpha}^{\delta} - \varphi^{\dagger}}_{\cH}$ 
by some rule for the choice of regularization
parameter has been studied and established well,
see \textbf{\cite{Engl96}, \cite{Isakov06}, \cite{Kirsch11}} for the details.
This convergence is also known as the {\em total error} and is defined by

\begin{displaymath}
E(\varphi_{\alpha}^{\delta} , \varphi^{\dagger}) := \norm{\varphi_{\alpha}^{\delta} - \varphi^{\dagger}}_{\cH}.
\end{displaymath}
From this presentation, one must expect from the numerical experiments that
the most reliable solution will be provided by the algorithm which gives
the least total error value during the iteration. Aside from the convergence
analysis in the pre-image space, we will also focus on the convergence
in the image space by analysing the discrepancy between 
$\cT\varphi_{\alpha}^{\delta}$ and the measured data $f^{\delta},$
{\em i.e.} $\norm{\cT\varphi_{\alpha}^{\delta} - f^{\delta}}_{\cL^{2}(\cZ)}.$
According to well-known Morozov's discrepancy principle (MDP),
one must define a rule for the choice of the regularization parameter
in a way such that the following, with some fixed 
$1 < \underline{\tau} \leq \overline{\tau} < \infty,$

\beq
\label{discrepancy_pr_definition}
\alpha(\delta , f^{\delta}) \in \{ \alpha > 0 \mbox{ } : \mbox{ }
\underline{\tau}\delta \leq \norm{\cT\varphi_{\alpha(\delta , f^{\delta})}^{\delta} - f^{\delta}}_{\cL^2(\cZ)} \leq \overline{\tau}\delta \} 
\mbox{ for all given } (\delta , f^{\delta}),
\eeq
must hold. Our tests do not involve any implementation
of the discrepancy principle. However, it is still in the expectations
of our tests that after some some number of iteration steps,
the convergence rate $\norm{\cT\varphi_{\alpha}^{\delta} - f^{\delta}}_{\cL^{2}(\cZ)}$
in the image space is expected to remain constant.

The updated reconstruction $\varphi_{\alpha}^{\nu}$ will be produced
by different gradient based algorithms, see Section \ref{quasi-Newton}
for the details. Thus, significant decay in the norm of the gradient
of the functional is expected, {\em i.e.} 
$ \Vert\nabla F(\varphi_{\alpha}^{\nu})\Vert \leq \Vert\nabla F(\varphi_{\alpha}^{\nu - 1})\Vert$
at each iteration step $\nu = 1, 2, \cdots .$


\subsection{The synthetic profile}
\label{the_synthetic_profile}

The atmospheric physical facts behind the refractivity profile
of humidity fields can be found in \textbf{\cite{Kleijer04, Perler11}}.
The vertical profile of the refractivity $\varphi$ can be approximated
by an exponential function, (cf. \textbf{\cite[Eq. (17)]{Perler11}}),
with the empirically determined scale height parameters $H_{1c}$ and $H_{2c},$

\bea
\label{vertical_profile}
\varphi^{\dagger}(h) = \frac{N_{0}}{2} \left( \exp\left\{ - \frac{h}{H_{1c}} \right\} + \exp\left\{ - \frac{h}{H_{2c}} \right\} \right).
\eea
Linear functions of $x$ and $y$ would introduce gradients along these axes.
Periodical variations are modelled to define horizontal profile,

\bea
\label{horizontal_profile}
\varphi^{\dagger}(x,y) = N_{0} + \frac{N_{x}x}{\Delta_{x}} + \frac{N_{y}y}{\Delta_{y}} +
N_{1} \sin{\left( \frac{2\pi\mu_{x}x}{ \Delta_x } \right)} + 
N_{2} \cos{\left( \frac{2\pi\mu_{y}y}{ \Delta_y } \right)},
\eea
where $\Delta_x = x_{\mathrm{max}} - x_{\mathrm{min}}$ 
and $\Delta_y = y_{\mathrm{max}} - y_{\mathrm{min}},$ 
$N_{1}$ and $N_{2}$ are the amplitudes of the periodic variations, 
$\mu_{x}$ and $\mu_{x}$ are the corresponding frequencies which are normalized
to the $x$ and $y$ intervals. Combining everything one gets a 
three dimensional refractivity field with number of parameters

\bea
\label{true_profile}
\varphi^{\dagger}(x,y,h) = \frac{N_{0}}{2} \left[ N_{0} + \frac{N_{x}x}{\Delta_{x}} + \frac{N_{y}y}{\Delta_{y}} +
N_{1} \sin{\left( \frac{2\pi\mu_{x}x}{ \Delta_x } \right)} + 
N_{2} \cos{\left( \frac{2\pi\mu_{y}y}{ \Delta_y } \right)} \right]
\nonumber\\
\left( \exp\left\{ - \frac{h}{H_{1c}} \right\} + \exp\left\{ - \frac{h}{H_{2c}} \right\} \right).
\eea
For the parameters defined as $\mu_{x} = 4,$ $\mu_{y} = 6,$ $N_{0} = 350,$
$H_{1c} = 1,$ $H_{2c} = 7,$ $N_x = 30,$ $N_y = 50,$ $N_1$ and $N_2$
can be chosen in a way $N_{0} - N_{1} - N_{2} \geq 200$ and 
$N_{0} + N_{1} + N_{2} \leq 400.$
Below in Figure \ref{3-D_ProfilesNumerics}, true and the noisy solutions
can be seen for the numerical experiments. 

\begin{figure}[ht]
\centering
\includegraphics[height=5in,width=6.5in,angle=0]{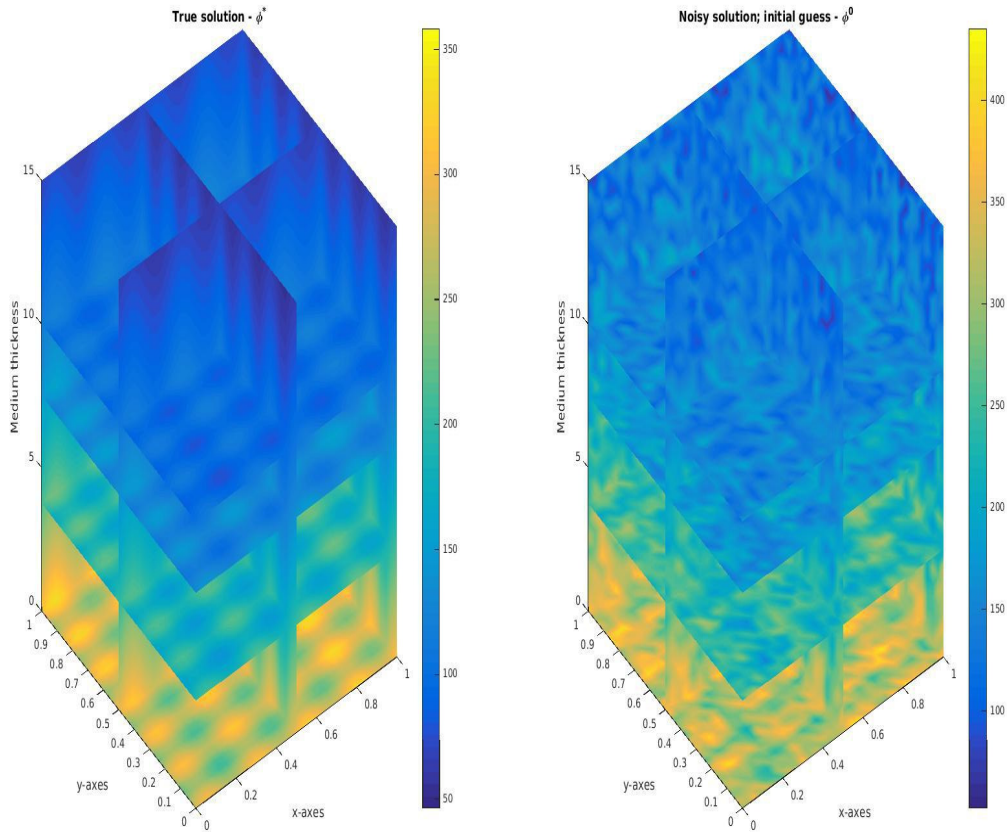}
\caption[Simulated true and noisy solutions for the numerical experiments.]
{\footnotesize Simulated true and noisy solutions for the numerical experiments. 
The domain $\Omega$ has been discretized by $30 \times 30 \times 30$ points.}
\label{3-D_ProfilesNumerics}
\end{figure}


\subsection{On the implementation of the forward operator}
\label{forward_op_imp}

\begin{figure}[ht]
\centering
\includegraphics[height=5in,width=6.5in,angle=0]{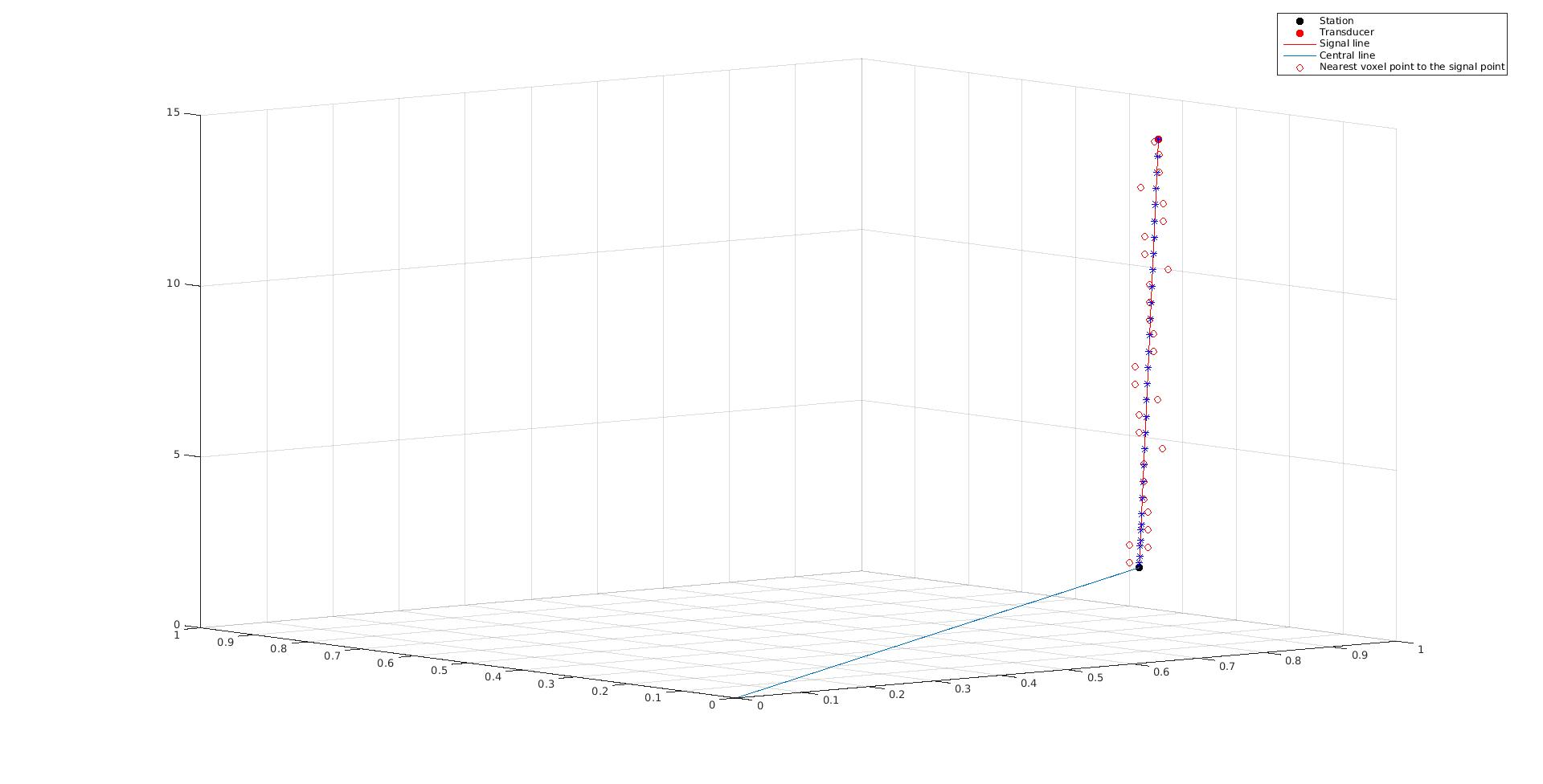}
\caption[Demonstration for the line integration procedure.]
{\footnotesize Demonstration for the line integration procedure.}
\label{line_integration}
\end{figure}

Thorough implementation and inversion of geodesic X-ray
transform has been studied in \textbf{\cite{Monard14}}.
Here, we focus on the linearized form of that regarding general
implementation. In the computerized environment, we are only capable of 
implementing discretized integration which has been introduced
in (\ref{disc_integration}).
In our implementation, this discretized integration is carried on according to 
{\em nearest neighbor search}, or {\em closest point search},
principle. To this end, discretization of each ray $\gamma$ 
is necessary. Owing to the parameter function 
$t(\epsilon) = \frac{\epsilon}{\sin(\tilde{\rho})},$
where $\epsilon \in [z_s , h^{\infty}],$
we are able to discretize $\gamma,$
see (\ref{signal_set}).  
For one ray, this discretization is illustrated in 
Figure \ref{line_integration} whereby blue stars denote the mesh 
points of the signal path $\gamma$ and 
the red circles are for the nearest points to the corresponding mesh 
point of $\gamma.$ Discretized line integration is carried 
on along those red circles. The implemented 
integration procedure seeks the nearest point 
to the corresponding interior point of 
$\gamma$ on the horizontal layer.
By the nearest point, we mean the closest
grid point of the area of interest $\Omega$ to the interior point
of the corresponding ray. This procedure can be described
mathematically as such; For any index $k \in \{ i , i + 1 \}$ where
$ i = 1, \cdots, N - 1,$ denote by $\vec{x}_{k}$ any grid point of 
our simulated area of interest $\Omega \subset \R^{3}.$ 
Interior point of any ray $\gamma_{[\mb{s} , \vec{\theta}]}$
is denoted by $\vec{\gamma}_{l}$ for $l = 1, \cdots, m,$
where $m$ is the number of the interior points.
Then, we seek the closest point $\vec{x}_{k} \in \Omega$
to the interior point $\vec{\gamma}_{l} \in \gamma_{[\mb{s},\vec{\theta}]}$ 
according to the finite dimensional maximum norm by 

\bea
\label{closest_distance}
d_{k}^{l} := \min_{i} \{ \Vert \vec{\gamma}_{l} - \vec{x}_{i} \Vert_{\infty} , \Vert \vec{\gamma}_{l} - \vec{x}_{i + 1} \Vert_{\infty} \} 
\mbox{, for } k \in \{ i , i + 1 \} \mbox{ where }  i = 1, \cdots, N - 1.
\eea
The pointwise density value at the corresponding point 
$\vec{x}_{k} \in \Omega \subset \R^{3}$ is 
$\varphi_{k} = \varphi(\vec{x}_{k}),$ for $k \in \{ i , i + 1 \}$
where $ i = 1, \cdots, N - 1.$ 
Eventually, the true measurement vector $f_{j}^{\dagger}$ 
for the corresponding ray is calculated by

\bea
[T\varphi]_{j} = \sum_{i = 1}^{N} \varphi_{i} w_{i,j} = f_{j}^{\dagger}.
\eea

%

\section{Numerical Results and Review on the Algorithms}
\label{numerics}

Since application of smooth TV is a new regularization
strategy for this particular problem, it is expected to obtain
some reasonable reconstruction. We will also realize usual
facts in regularization theory. Firstly, this problem can also
be interpreted as another sparse reconstruction. Therefore,
measurement number (number of the signal) will impact
on the convergence rate in the pre-image space. We will
demonstrate this by the relative error value of the reconstruction.
Secondly, as well known by the usual regularization theory
for the inverse ill-posed problems \textbf{\cite{Engl96}},
noise amount defined in the image space will also have
impact on the convergence rate in the pre-image space.
This latter case will also be demonstrated by visualizing
the relative error value of the reconstruction.


\subsection{Quasi-Newton Methods}
\label{quasi-Newton}

Much of the technical and scientific details of this section
can be found in \textbf{\cite{Luke01, NocedalWright99}}.
With a positive definite symmetric approximate Hessian
$H_{\alpha}^{\nu}$ and properly chosen step-length parameter
$\eta,$ one obtains the update $\varphi_{\alpha}^{\nu + 1}$ 
by a quasi-Newton method as such,

\bea
\label{quasi-Newton}
\varphi_{\alpha}^{\nu + 1} = \varphi_{\alpha}^{\nu} - \eta H_{\alpha}^{\nu} \nabla F_{\alpha}^{\nu} .
\eea
Here, $F_{\alpha}^{\nu} = F_{\alpha}(\varphi_{\alpha}^{\nu} , f^{\delta}).$
Algorithms that provide the approximate solution $\varphi_{\alpha}^{\nu + 1}$
are called {\em gradient based algorithms}, 
(cf. \textbf{\cite{BeckTeboulle09}} and \textbf{\cite[Eq 2.1]{BeckTeboulle09_2}}).


\subsection{Lagged Diffusivitiy Fixed Point Iteration - (LDFP)}

The favourite regularization strategy of this work is TV
regularization. Therefore, we would like to begin with
one of the simplest algorithms to illustrate our regularized
solution. LDFP, \textbf{\cite{Vogel02, VogelOman96}}, 
is also in the class of quasi-Newton search direction
algorithm. Since the Fr\'{e}chet differentiable functional 
$F_{\alpha}(\varphi , f^{\delta})$ is defined by

\begin{displaymath}
F_{\alpha}(\varphi , f^{\delta}) := \frac{1}{2} \norm{\cT\varphi - f^{\delta}}_2^2 
+ \alpha \int_{\Omega} \sqrt{\vert \nabla\varphi(x) \vert^2 + \beta} dx,
\end{displaymath}
then LDFP is given by the following scheme,

\bea
\label{LDFP_scheme}
\varphi_{\alpha}^{\nu + 1} & = & \varphi_{\alpha}^{\nu} + 
(\cT^{\ast}\cT + \alpha L(\varphi_{\alpha}^{\nu}))^{-1} \nabla F_{\alpha}(\varphi , f^{\delta})
\nonumber\\
& = & (\cT^{\ast}\cT + \alpha L(\varphi_{\alpha}^{\nu}))^{-1} \cT^{\ast}f^{\delta} = 
\cR_{\alpha}(\varphi_{\alpha}^{\nu}), \hspace{.2cm} \nu = 0, 1, \cdots ,
\eea
where,

\begin{displaymath}
L(\varphi_{\alpha}^{\nu}) := - \nabla^{\ast} \cdot \left( \frac{\nabla}
{( \vert \nabla \varphi_{\alpha}^{\nu} \vert^2 + \beta)^{1/2}} \right) .
\end{displaymath}
Comparison between (\ref{LDFP_scheme}) and (\ref{quasi-Newton})
yields that in the LDFP scheme the step-length $\eta = 1,$ and
the approximate Hessian is defined by

\bea
\label{Hessian_LDFP}
H_{\alpha}^{\nu} := \cT^{\ast}\cT + \alpha L(\varphi_{\alpha}^{\nu}).
\eea
Direct implementation of the scheme (\ref{LDFP_scheme}) would still be a costly
iteration procedure since $L(\varphi)$ is 
highly nonlinear. Then, according to \textbf{\cite[Algorithm 8.2.3]{Vogel02}}, 
the update $\varphi_{\alpha}^{\nu + 1}$
is produced after the following linearization steps;

\mbox{ }

\underline{\textbf{LDFP algorithm with smooth-TV penalty:}}

\mbox{ }

\begin{figure}{ }
  \center
  \includegraphics[height=5in,width=6.5in,angle=0]{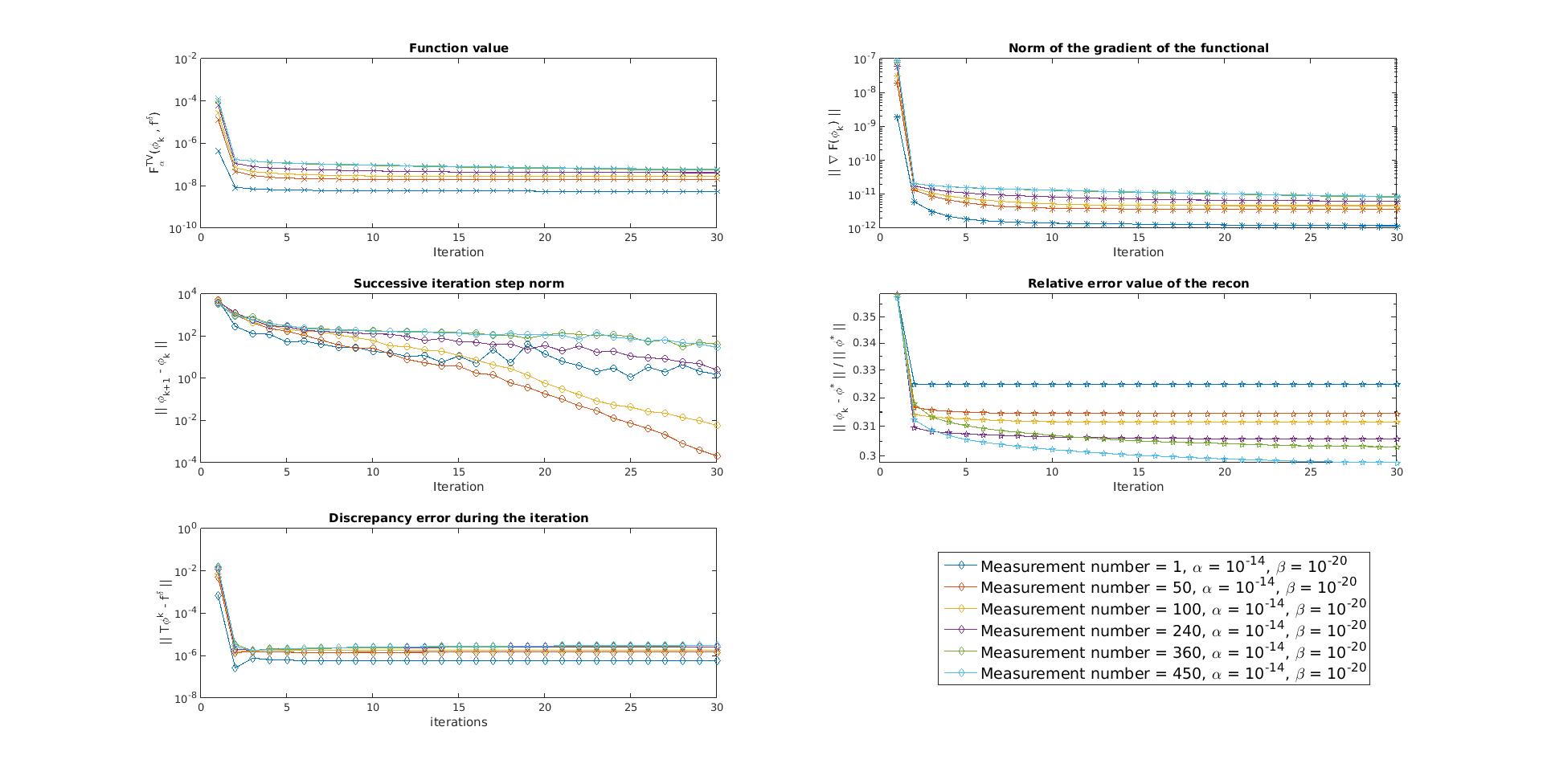}
  \caption[LDFP numerical convergence results per different number of the measurements.]
  {\footnotesize LDFP numerical convergence results per different number of the measurements,
  $\{1, 50, 100, 240, 360, 450\}.$ Regularization parameter is chosen according to the stable behaviour of the discrepancy 
  $\norm{\cT\varphi_{\alpha}^{\nu} - f^{\delta}}$ after each iteration step $\nu = 1, 2, 3, \cdots .$ }
  \label{LDFP_TV_convergence_benchmark1}
\end{figure}

\begin{figure}{ }
  \center
  \includegraphics[height=5in,width=6.5in,angle=0]{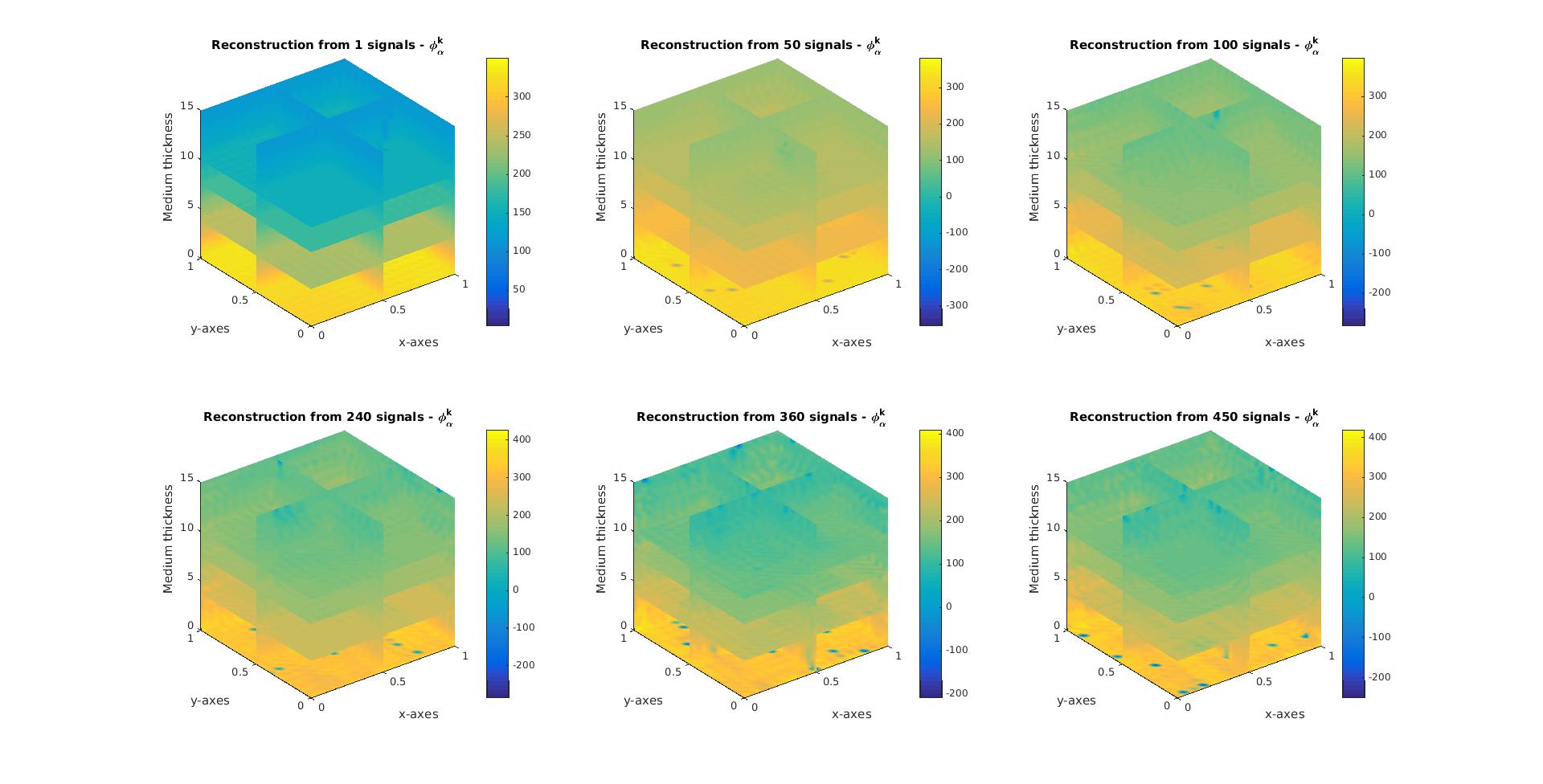}
  \caption[LDFP numerical reconstruction results per different number of the measurements.]
  {\footnotesize LDFP numerical reconstruction results per different number of the measurements,
  $\{1, 50, 100, 240, 360, 450\}.$ }
  \label{LDFP_TV_recon_benchmark1}
\end{figure}

\begin{algorithmic}

\State 1. Compute $L^{\nu} := L(\varphi_{\alpha}^{\nu})$ anisotropic Laplacian;

\State 2. Compute $g^{\nu} := \mathcal{T}^{\ast}(\mathcal{T}\varphi_{\alpha}^{\nu} - f^{\delta}) + 
\alpha L^{\nu}\varphi_{\alpha}^{\nu}$ gradient step;

\State 3. Compute $H_{\alpha}^{\nu} := \mathcal{T}^{\ast}\mathcal{T} + \alpha L^{\nu}$ 
approximate Hessian;

\State 4. Solve $H_{\alpha}^{\nu} s^{\nu + 1} = -g^{\nu}$ quasi-Newton step;

\State 5. Update $\varphi_{\alpha}^{\nu + 1} = \varphi_{\alpha}^{\nu} + s^{\nu + 1};$

\end{algorithmic}
In our experiments, we use usual CGNE for solving the inner system 
$H_{\alpha}^{\nu} s^{\nu} = g^{\nu},$ see \textbf{\cite{Hanke95}}.
In the Figure \ref{LDFP_TV_convergence_benchmark1}, we present the numerical results of 
LDFP algorithm per different number of the measurements. 
We run the algorithm only for $30$ iteration steps to understand its behaviour.
Reconstructions that are the results of LDFP algorithm per different number
of measurements are presented in Figure \ref{LDFP_TV_recon_benchmark1}.

\subsection{Quasi-Newton method for large-scale problems}
\label{large_scale_quasi_newton}

The quasi-Newton methods cannot be directly applicable to large
optimization problems because their approximations to the Hessian
or its inverse are usually dense. The storage and computational
requirements grow in proportion to $N^2,$ and become excessive
for large $N.$ In order to overcome this difficulty, limited-memory
quasi-Newton methods have been introduced, 
\textbf{\cite{Luke01, NocedalWright99}}. Here, we particularly
focus on {\em limited memory BFGS} (L-BFGS) algorithm.

By applying a quasi-Newton method, finding the optimum solution to the
minimization problem (\ref{problem}), amounts to solving
secant equation given by

\bea
\label{secant_eq}
B^{\nu + 1} s_{\nu} = y_{\nu},
\eea
where

\bea
s_{\nu} = \varphi_{\alpha}^{\nu + 1} - \varphi_{\alpha}^{\nu}, \mbox{ } y_{\nu} = \nabla F_{\alpha}^{\nu + 1} - \nabla F_{\alpha}^{\nu} .
\eea
Here, $F_{\alpha}^{\nu} = F_{\alpha}(\varphi_{\alpha}^{\nu} , f^{\delta}).$
In (\ref{secant_eq}), the matrix $B^{\nu + 1}$ is a positive definite symmetric
approximation to the true Hessian of the cost functional $F_{\alpha}.$

\bea
\varphi_{\alpha}^{\nu + 1} = \varphi_{\alpha}^{\nu} - \eta^{\nu} H^{\nu} \nabla F_{\alpha}^{\nu} .
\eea
Here the aproximate Hessian $H^{\nu}$ is updated by ,

\bea
H^{\nu + 1} = V_{\nu}^{T} H^{\nu} V_{\nu} + \rho_{\nu} s_{\nu} s_{\nu}^{T}
\eea
with,
\bea
\rho_{\nu} = \frac{1}{y_{\nu}^{T}s_{\nu}}, \mbox{ } V_{\nu} = I - \rho_{\nu} y_{\nu} s_{\nu}^{T}.
\eea


\subsection{L-BFGS with trust region}
\label{L_BFGS_trust_regions}

Robustness of L-BFGS algorithm is provided by trust regions, 
\textbf{\cite[p. 232]{DennisSchnabel83}\cite[p. 91]{Luke01}}.
The trust region is the set of all points, \textbf{\cite{ConnGouldToint00}},

\begin{displaymath}
\cB^{\nu} := \{\varphi_{\alpha}^{\nu} \mbox{ } \vert \mbox{ } \norm{\varphi - \varphi_{\alpha}^{\nu}} \leq \Delta^{\nu}\}.
\end{displaymath}
Trust-region subproblem with trust-region radius $\Delta^{\nu}$
has been described well in \textbf{\cite[p. 94]{Luke01}}.

\begin{figure}[ht]
\center
\includegraphics[height=4.5in,width=6.5in,angle=0]{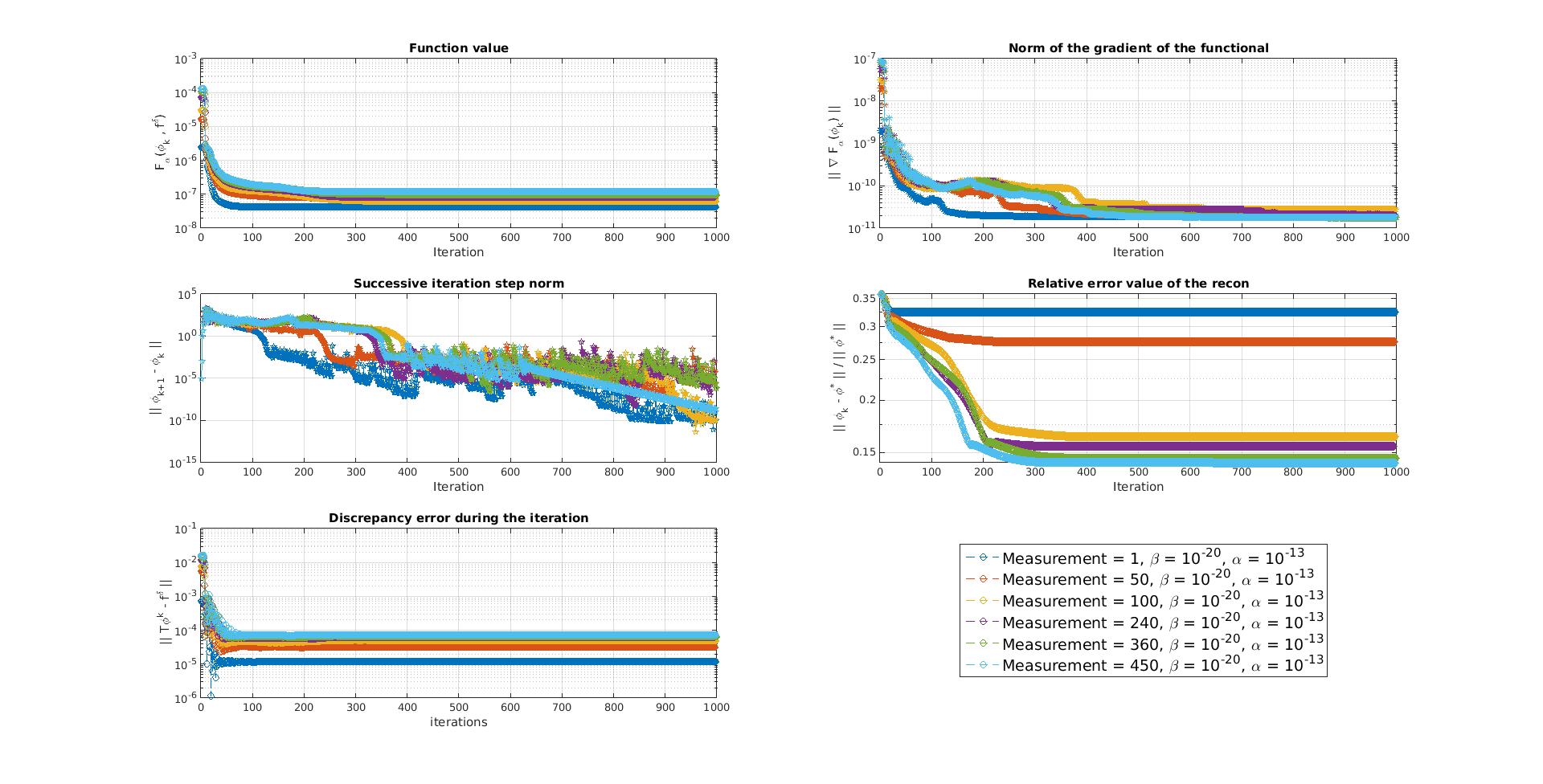}
  \caption[SAMSARA with TV gradient step numerical convergence results per measurement.]
  {\footnotesize SAMSARA with TV gradient step numerical results per measurement. We have conducted our experiment
  in the software SAMSARA for the measurement number $\{1, 50, 100, 240, 360, 450 \}.$}
\label{SAMSARA_TV_measurement_benchmark1}
\end{figure}

\begin{figure}[ht]
\center
\includegraphics[height=4.5in,width=6.5in,angle=0]{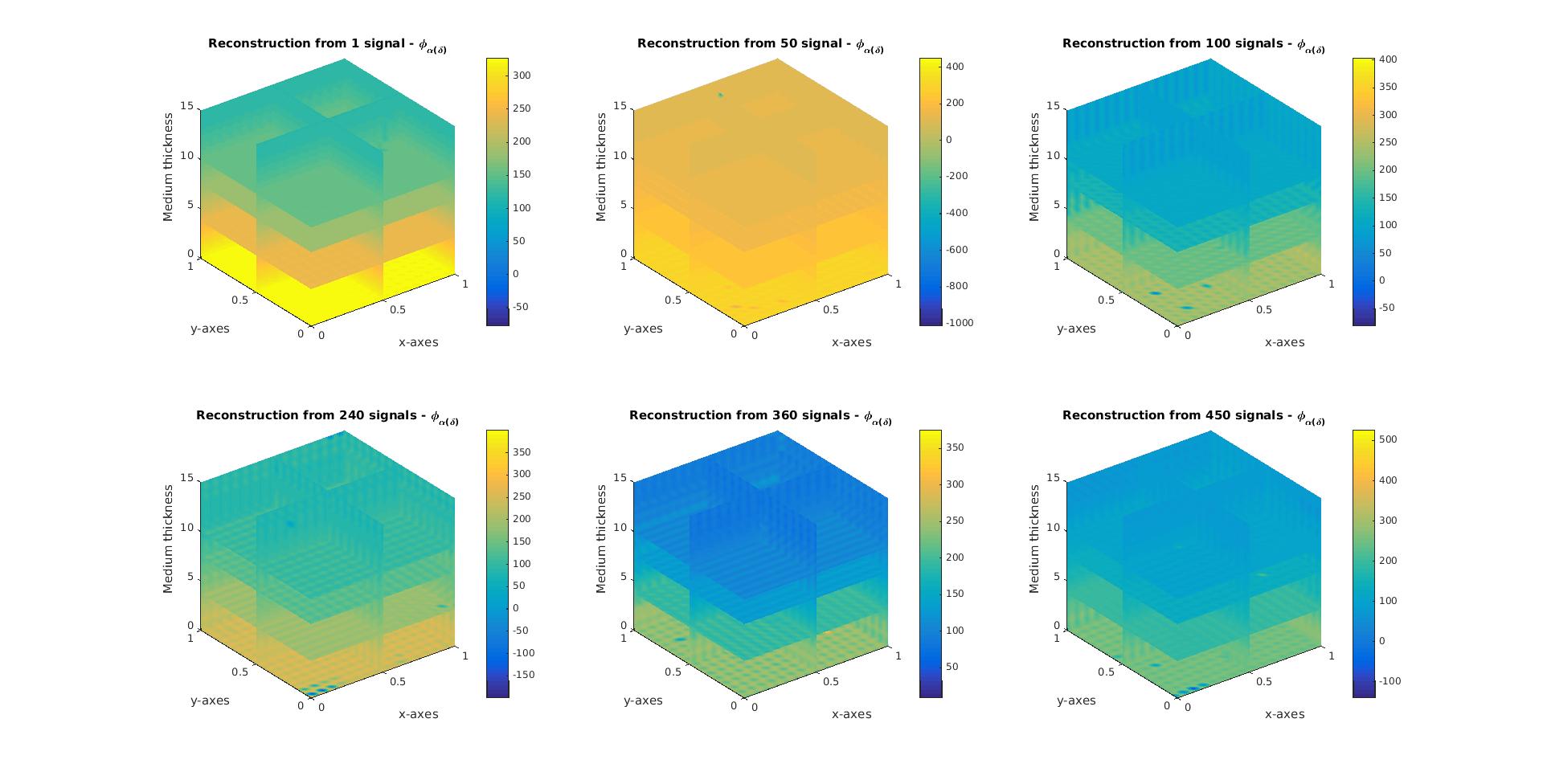}
  \caption[SAMSARA numerical reconstruction results per measurement.]
  {\footnotesize SAMSARA with TV gradient step numerical reconstruction results per measurement 
  $\{1, 50, 100, 240, 360, 450 \}.$ Fixed regularization parameter $\alpha = 10^{-13}$ has been determined
  according to the behaviour in the discrepancy.}
\label{SAMSARA_TV_measurement_benchmark2}
\end{figure}

\begin{figure}[ht]
\center
\includegraphics[height=4.5in,width=6.5in,angle=0]{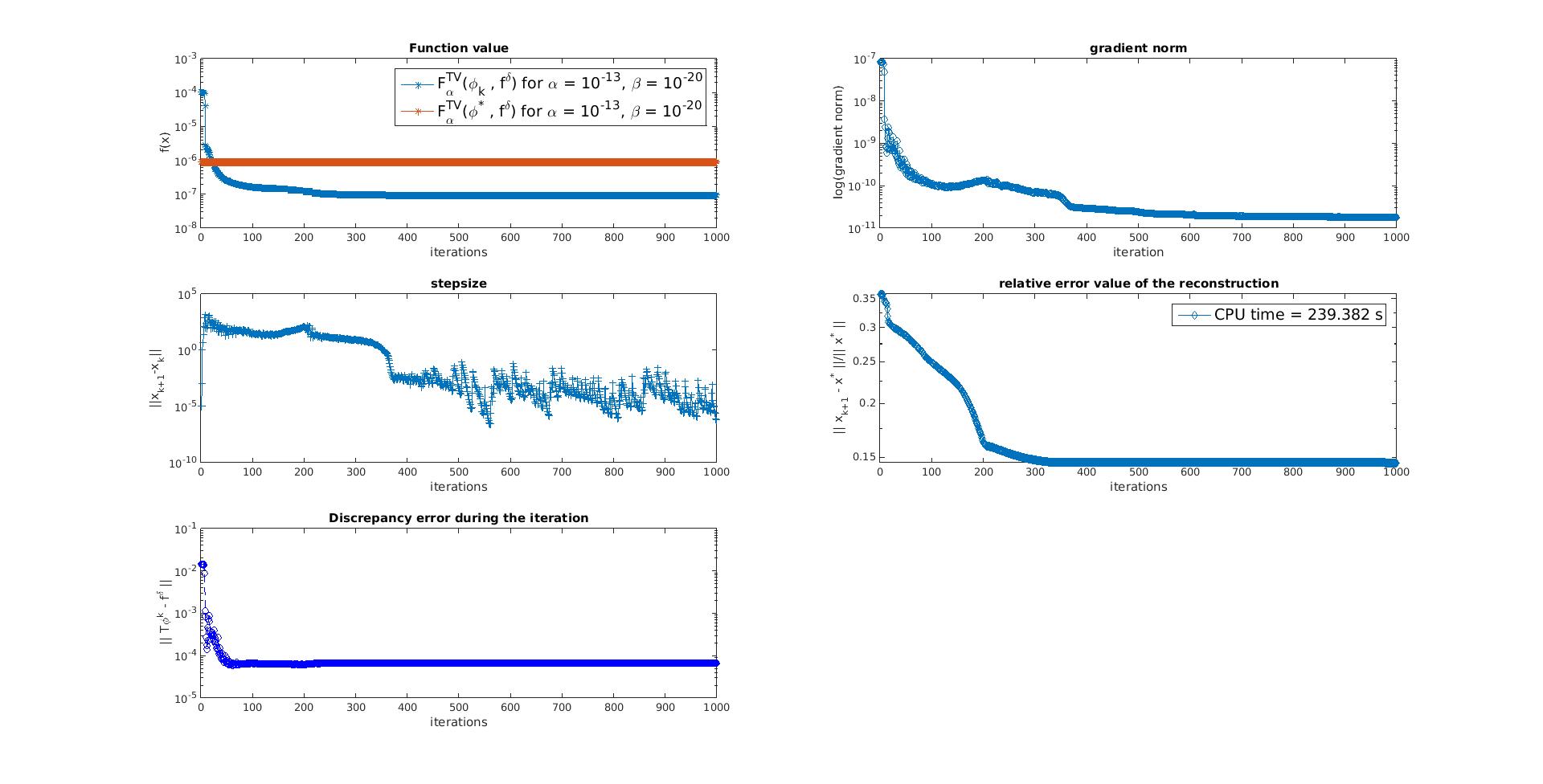}
  \caption[SAMSARA numerical convergence from 360 signals.]
  {\footnotesize SAMSARA numerical convergence from 360 signals.}
\label{SAMSARA_TV_360opt}
\end{figure}

\begin{figure}[ht]
\center
\includegraphics[height=4.5in,width=6.5in,angle=0]{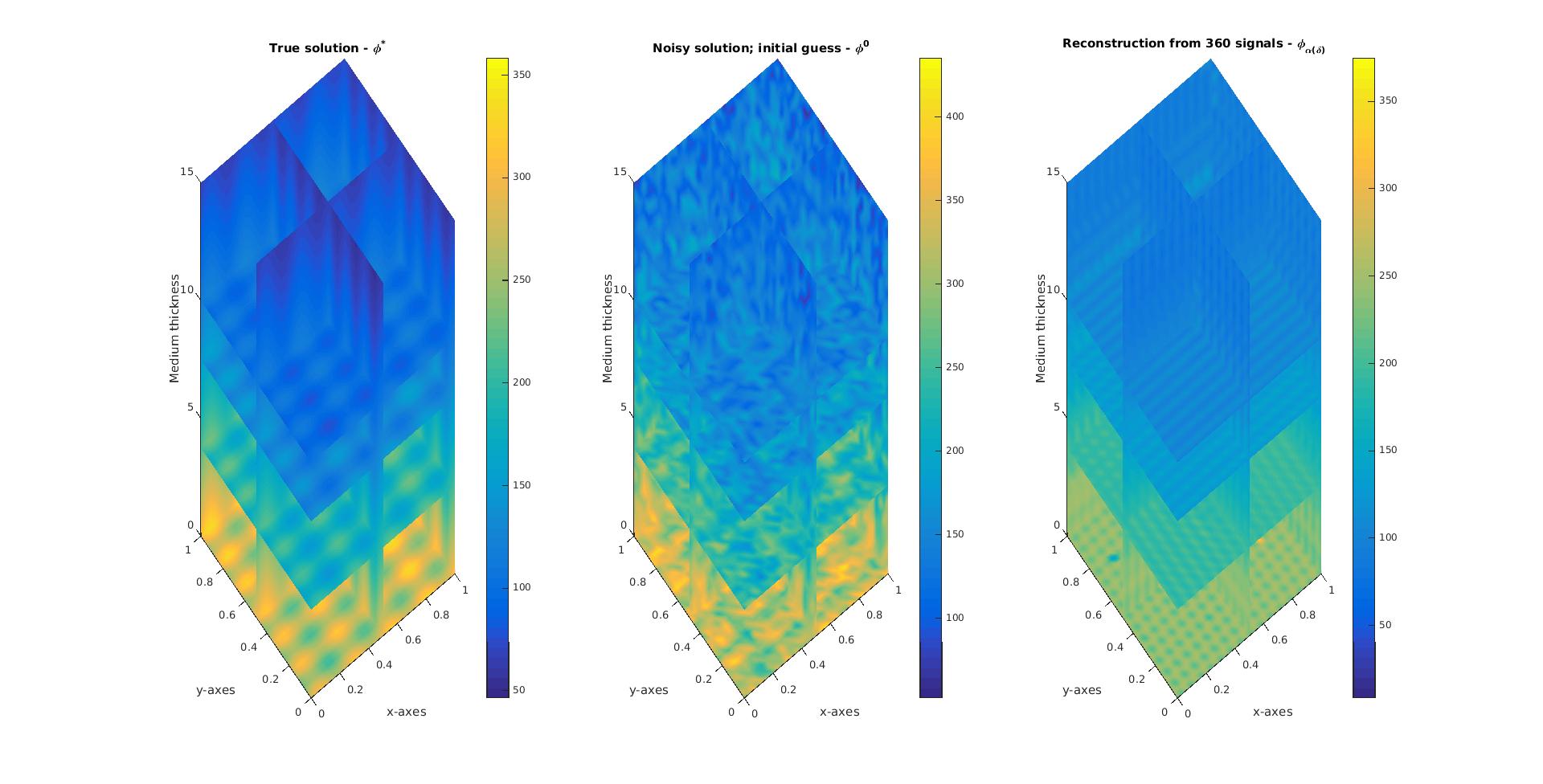}
  \caption[SAMSARA numerical reconstruction from 360 signals.]
  {\footnotesize SAMSARA numerical reconstruction from 360 signals.}
\label{SAMSARA_TV_recon1}
\end{figure}

\begin{figure}[ht]
\center
\includegraphics[height=4.5in,width=6.3in,angle=0]{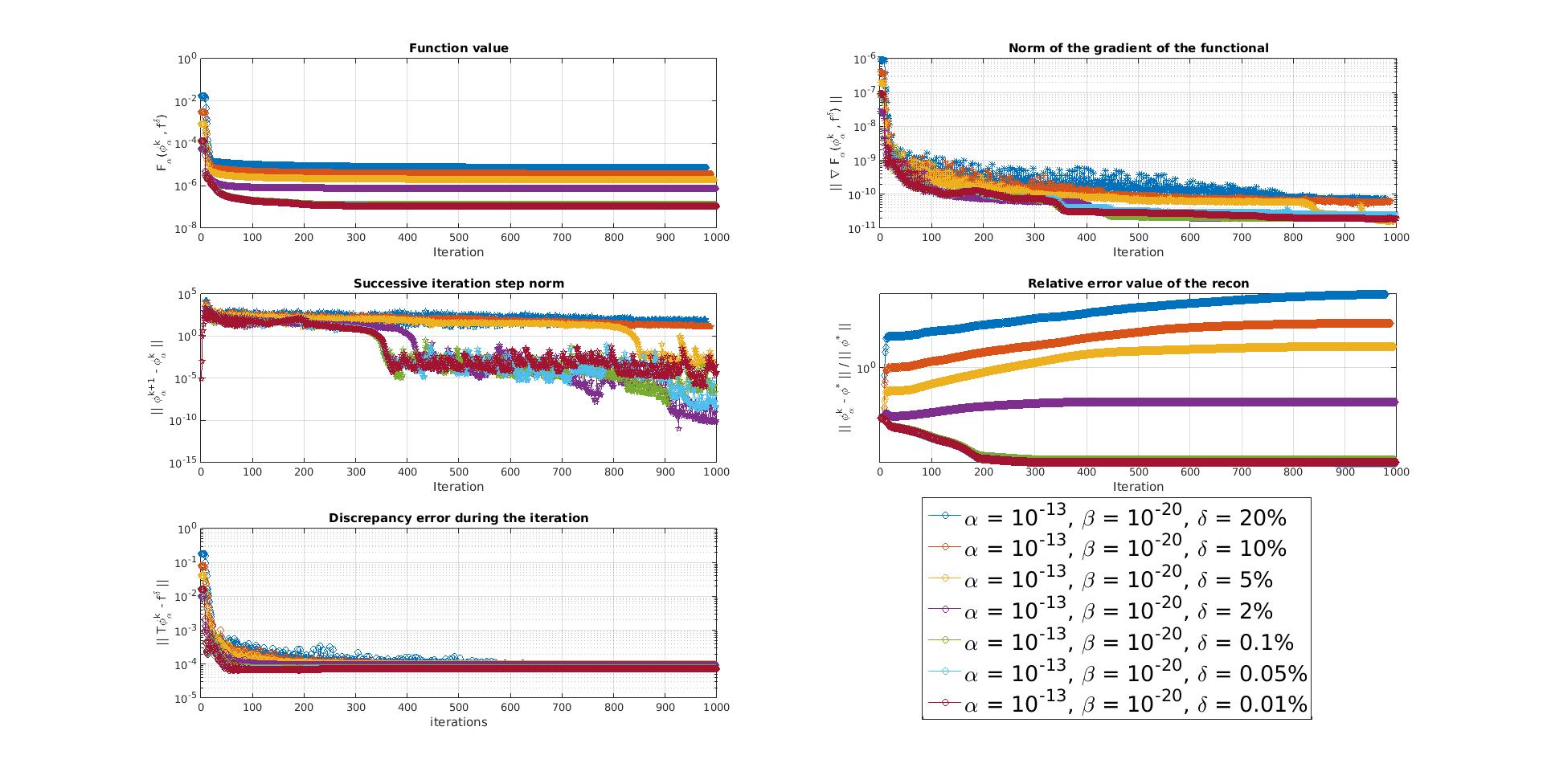}
  \caption[SAMSARA with TV gradient step convergence analysis per different noise amount.]
  {\footnotesize SAMSARA with TV gradient step convergence analysis per different noise amount 
  $\{ 20 \%, 10 \%, 5 \%, 2 \%, 0.1 \%, 0.005 \%, 0.001 \% \}.$ Convergence in the pre-image space
  begins with $\delta \leq 2 \%.$}
\label{SAMSARA_TV_noise_conv}
\end{figure}

\begin{figure}[ht]
\center
\includegraphics[height=4.5in,width=6.3in,angle=0]{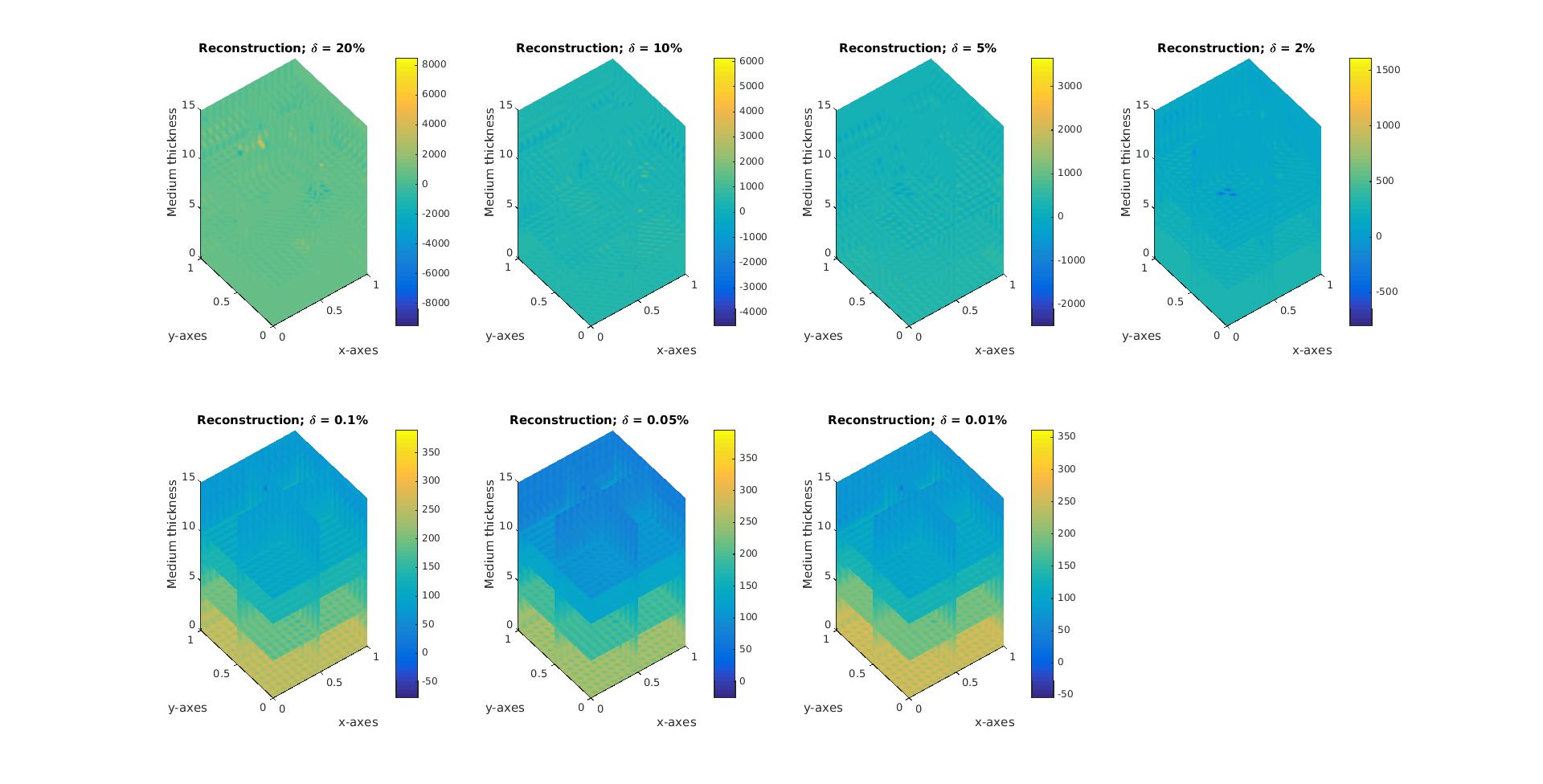}
  \caption[SAMSARA with TV gradient step different reconstruction per different noise amount.]
  {\footnotesize SAMSARA with TV gradient step different reconstruction per different noise amount.
  Corresponding numerical convergence analysis has been provided in Figure \ref{SAMSARA_TV_noise_conv}.}
\label{SAMSARA_TV_noise_recon}
\end{figure}

We provide the optimized solution for our problem (\ref{problem})
from  {\em trust region L-BFGS} algorithm by employing a novel reverse-communication large-scale nonlinear 
optimization software SAMSARA, \textbf{\cite{Luke}, \cite[Subsection 5.2.3]{Luke01}}. 
We demonstrate different solution per different measurement number, $\{1, 50, 100, 240, 360, 450 \},$
in the figures \ref{SAMSARA_TV_measurement_benchmark1}, \ref{SAMSARA_TV_measurement_benchmark2}, 
\ref{SAMSARA_TV_360opt}, \ref{SAMSARA_TV_recon1}. It is observed
better and more stable convergence rate in the pre-image space with
the more measurement number in the image space. Furthermore,
the figures \ref{SAMSARA_TV_noise_conv} and \ref{SAMSARA_TV_noise_recon} demonstrate
convergence in the pre-image/image spaces with varying amount of noise, 
$\delta \in \{ 20 \%, 10 \%, 5 \%, 2 \%, 0.1 \%, 0.005 \%, 0.001 \% \}.$
As a common expectation from an inverse ill-posed problem,
the less amount of noise in the image space provides better and stable
convergence rates in the pre-image space.

\subsection{Further reconstructions by L-BFGS with quadratic penalty Term}
\label{further_recon}

\begin{figure}{ }
  \center
  \includegraphics[height=4.5in,width=6.5in,angle=0]{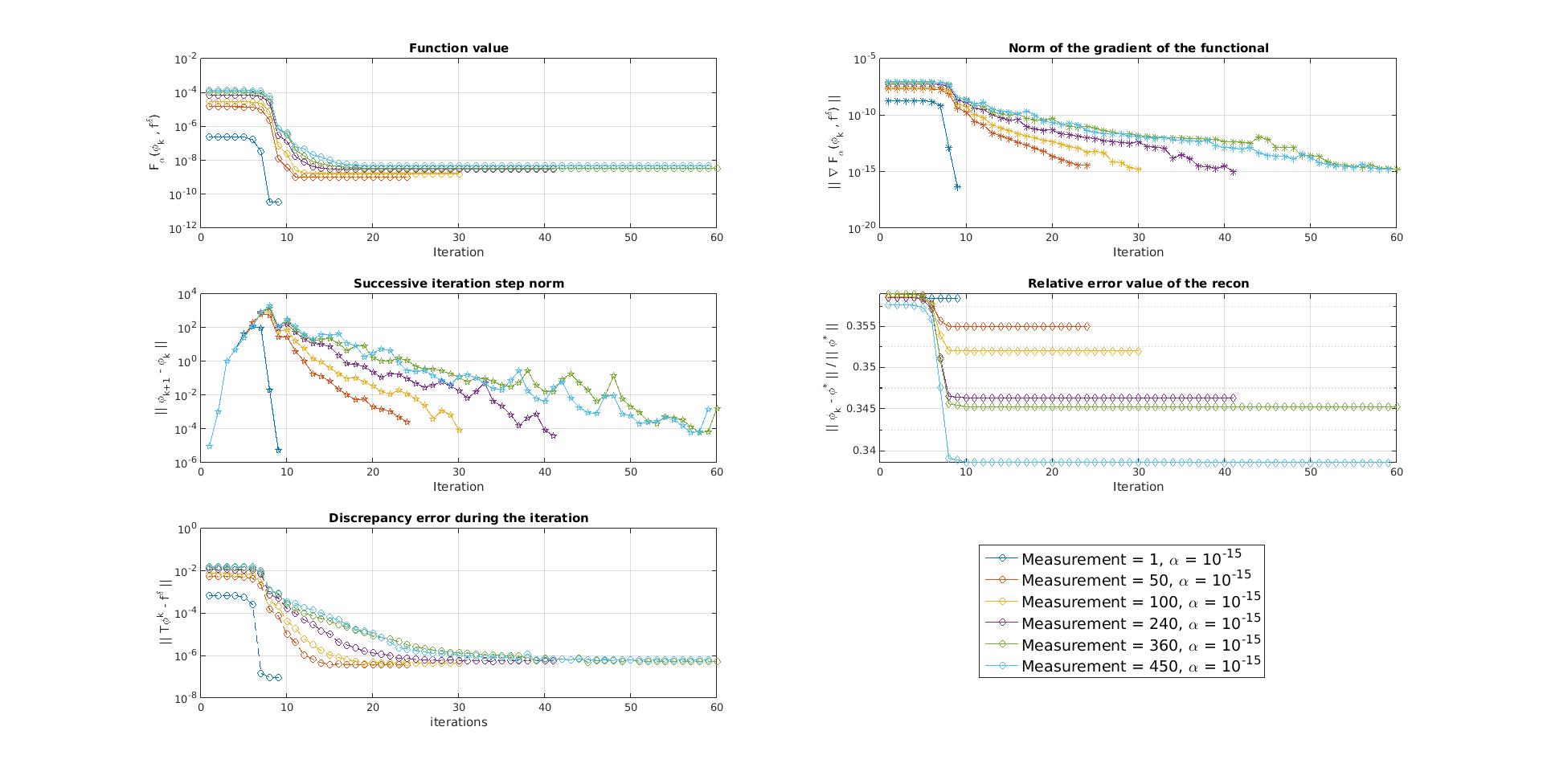}
  \caption[SAMSARA with Tikhonov gradient step numerical convergence results per measurement.]
  {\footnotesize SAMSARA results from quadratic Tikhonov gradient step numerical results per measurement. 
  We have conducted our experiment in the software SAMSARA for the measurement number $\{1, 50, 100, 240, 360, 450 \}.$ 
  Fixed regularization parameter $\alpha = 10^{-15}$ has been determined
  according to the behaviour in the discrepancy after each iteration step $\nu = 1, 2, \cdots $. }
  \label{SAMSARA_Tkh_measurement_benchmark1}
\end{figure}

\begin{figure}[ht]
\center
\includegraphics[height=4.5in,width=6.5in,angle=0]{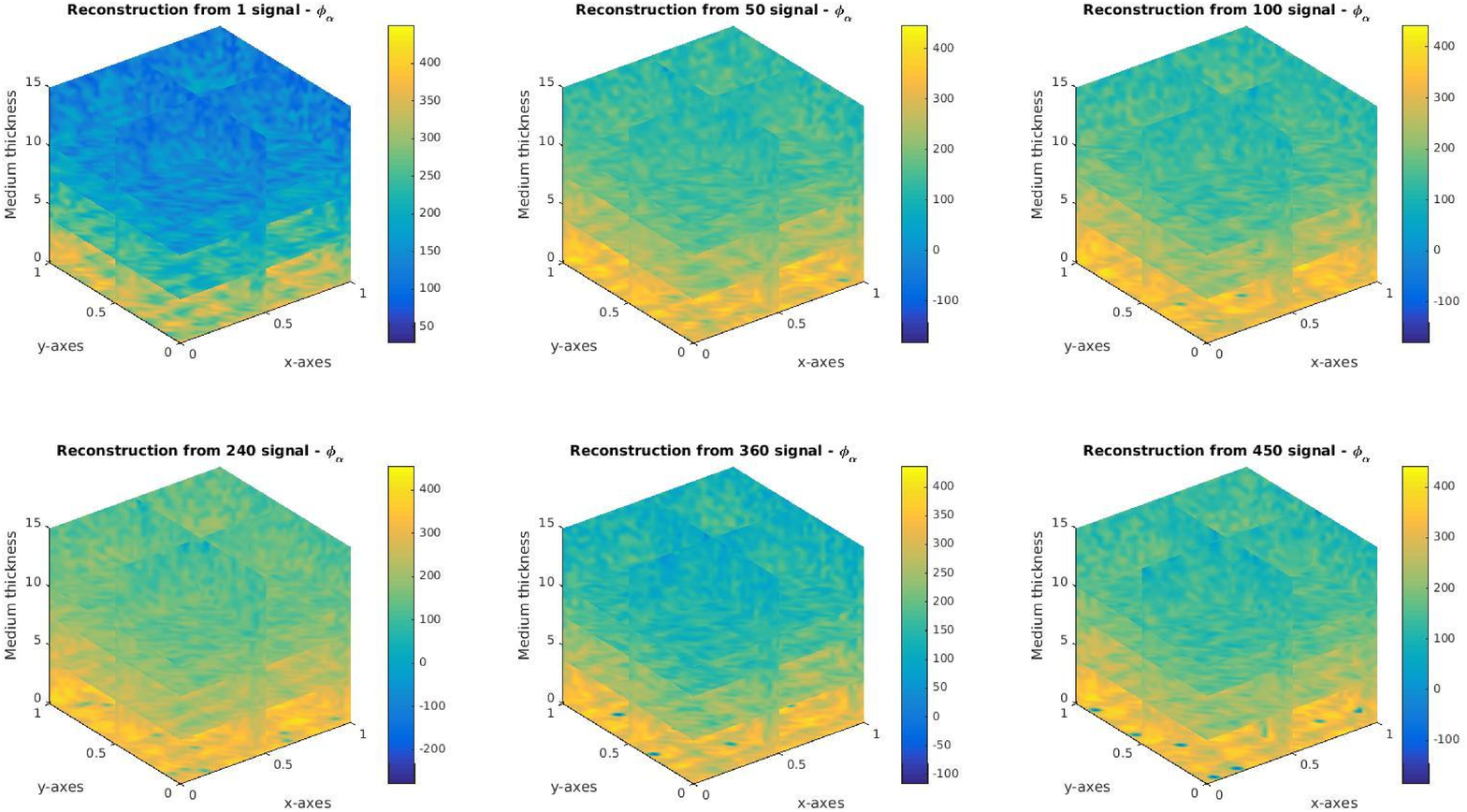}
  \caption[SAMSARA with quadratic penalty term gradient step numerical reconstruction results per measurement.]
  {\footnotesize SAMSARA result from quadratic Tikhonov gradient step numerical reconstruction results per measurement 
  $\{1, 50, 100, 240, 360, 450 \}$ for the fixed regularization parameter $\alpha = 10^{-15}$ }
\label{SAMSARA_Tkh_measurement_benchmark2}
\end{figure}

As a concrete demonstration for the effectiveness of TV type reconstruction,
we also seek approximate minimizer for the classical Tikhonov type functional 
with the quadratic term below,

\bea
\label{quadratic_cost_functional}
F_{\alpha}(\varphi , f^{\delta}) = \frac{1}{2}\norm{\cT\varphi - f^{\delta}}_{\cL^{2}(\cZ)}^2 + 
\alpha \frac{1}{2}\norm{\varphi - \varphi^{(0)}}_{\cL^{2}(\Omega)}^2 ,
\eea
with some given initial guess $\varphi^{0}.$ 
Here, we only present numerical results produced by SAMSARA, \textbf{\cite{Luke}},
with quadratic Tikhonov type objective functional.
Then the gradient step of the objective functional in (\ref{quadratic_cost_functional}) 
to be implemented is

\beq
\nabla F_{\alpha}(\varphi , f^{\delta}) = \cT^{\ast}(\cT\varphi - f^{\delta}) + \alpha(\varphi - \varphi^{(0)}) .
\eeq
We again run our tests with different number of measurements $\{1, 50, 100, 240, 360, 450 \}$
with sufficiently small amount of noise $\delta.$
Numerical convergence for each reconstruction is presented in
Figure \ref{SAMSARA_Tkh_measurement_benchmark1}. Each reconstruction is presented
in Figure \ref{SAMSARA_Tkh_measurement_benchmark2}. 


\section{Benchmark: LDFP vs SAMSARA With Smoothed-TV Penalty}
\label{LDFP_VS_SAMSARA}

\begin{figure}{ }
  \includegraphics[height=4.5in,width=6.5in,angle=0]{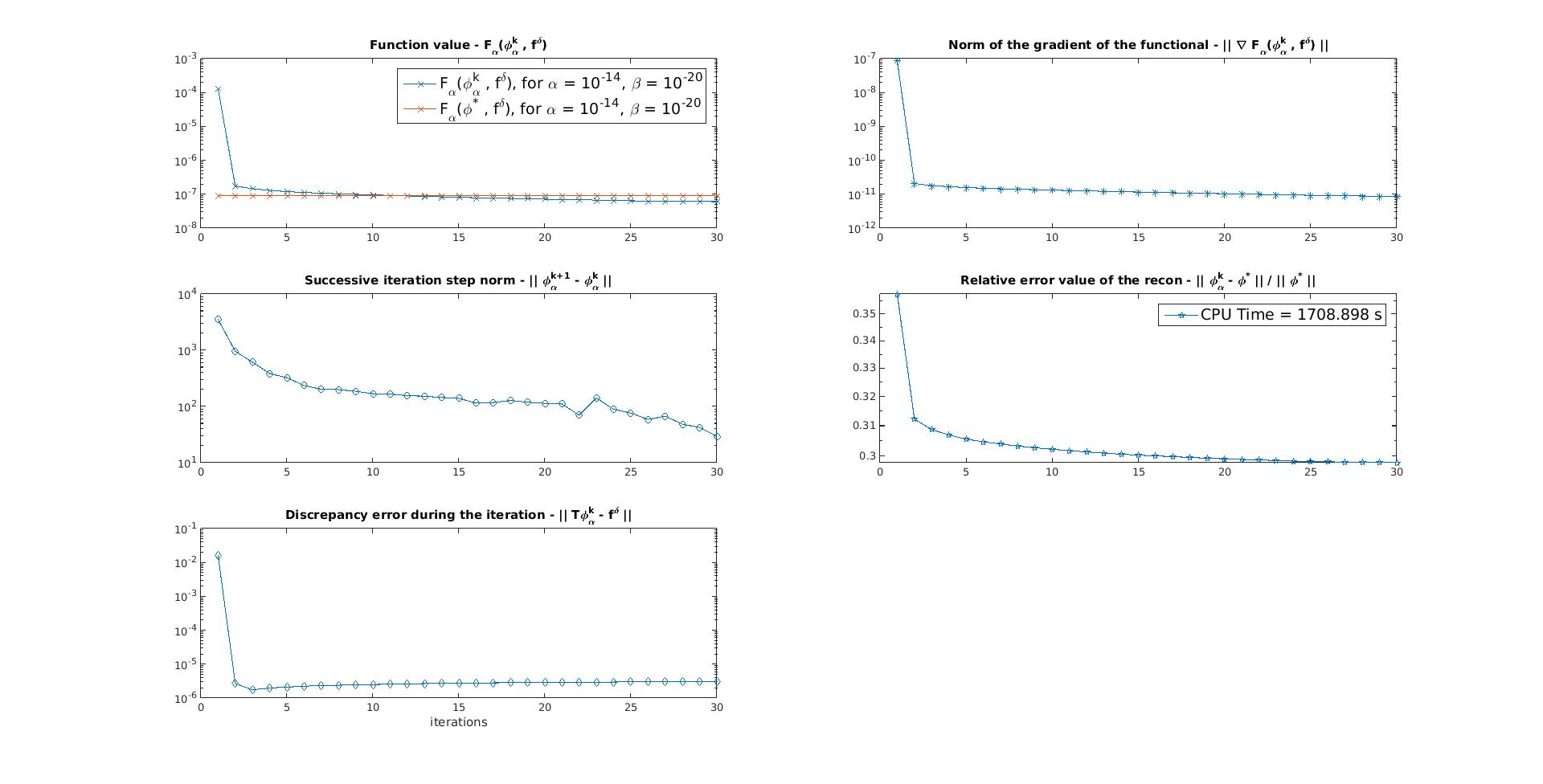}
  \caption[LDFP numerics from 450 measurements.]
  {\footnotesize LDFP numerics from 450 measurements.}
  \label{LDFP_numerics_from_450}
\end{figure}

\begin{figure}{ }
  \includegraphics[height=4.5in,width=6.5in,angle=0]{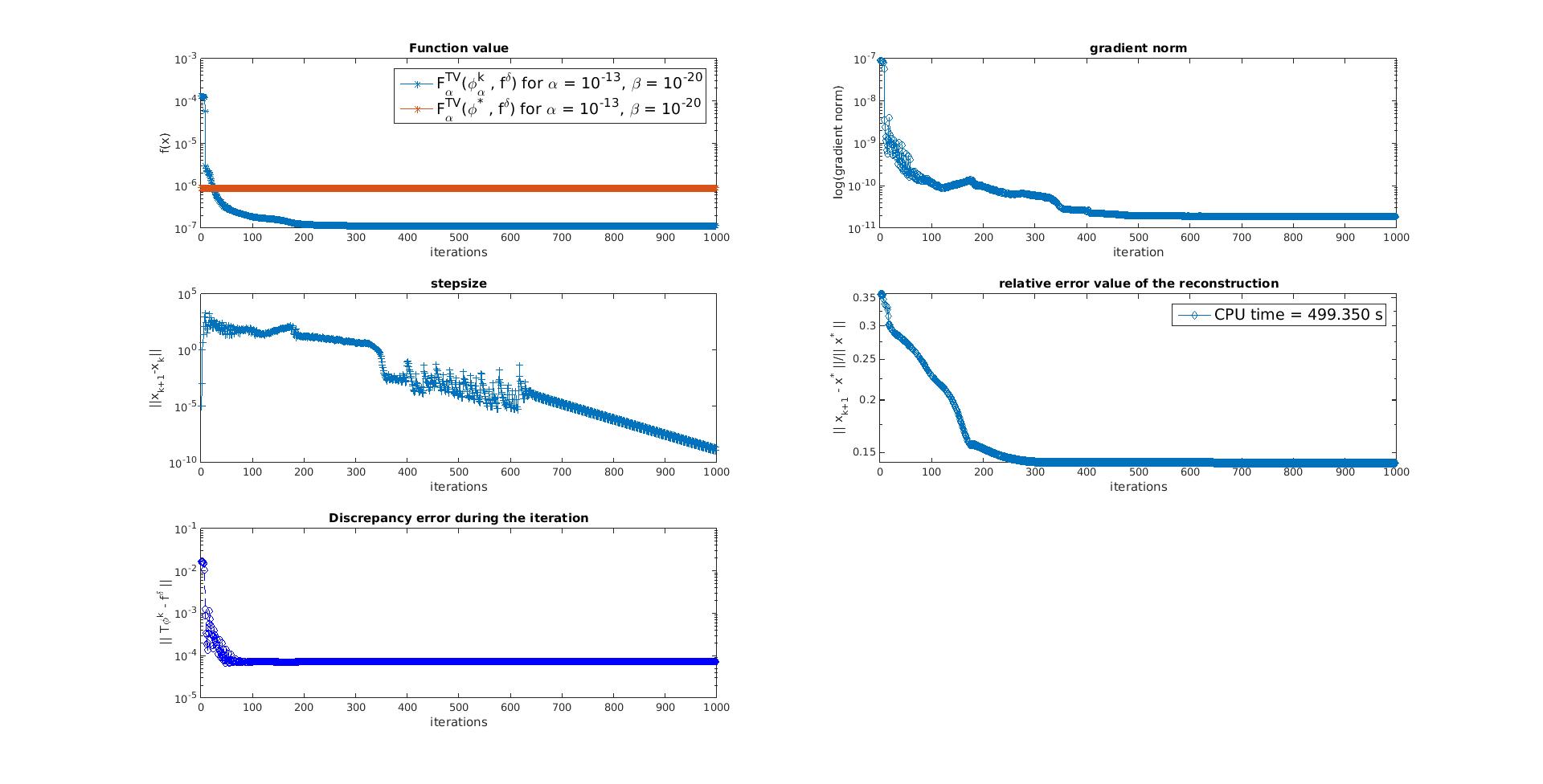}
  \caption[SAMSARA numerics from 450 measurements.]
  {\footnotesize SAMSARA with smoothed-TV gradient step numerics from 450 measurements.}
  \label{SAMSARA_numerics_from_450}
\end{figure}

\begin{figure}{ }
  \includegraphics[height=4.5in,width=6.5in,angle=0]{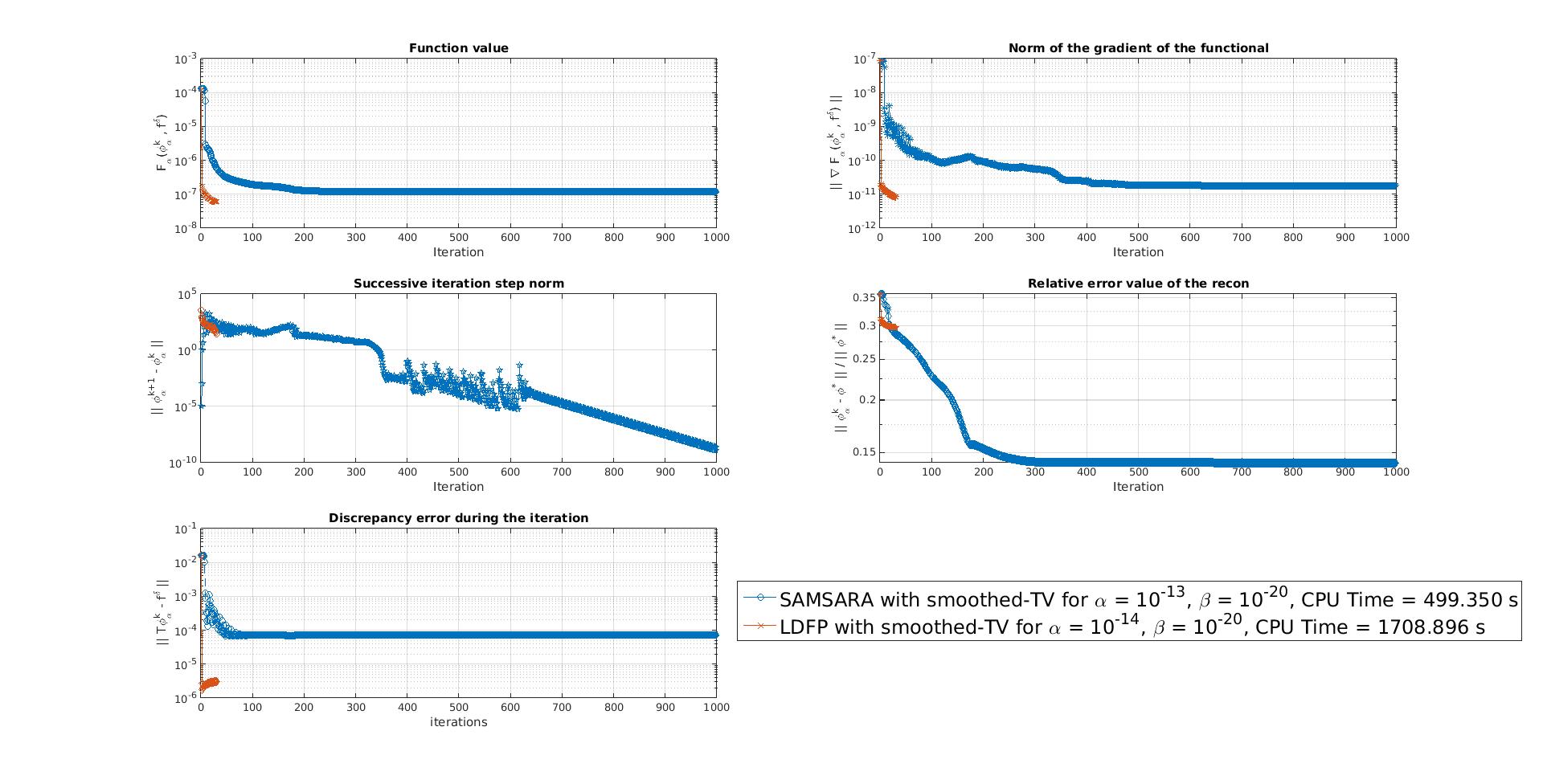}
  \caption[Benchmark: SAMSARA and LDFP numerics from 450 measurements.]
  {\footnotesize Benchmark: SAMSARA and LDFP numerics from 450 measurements. 
  Recall from the figures \ref{LDFP_numerics_from_450} and \ref{SAMSARA_numerics_from_450} that LDFP got executed
  only for 30 iteration steps whereas SAMSARA 1000 iteration steps. Comparing the CPU times for each tests, it would
  take SAMSARA only 14.9805 seconds to iterate 30 times. In other words, SAMSARA works 33 times faster than usual LDFP
  algorithm with smoothed-TV gradient step.}
  \label{SAMSARAvsLDFP_convergence_from_450}
\end{figure}

A CPU time based benchmark test between SAMSARA and LDFP both
associated with the smoothed-TV gradient step has been conducted,
see the figures \ref{LDFP_numerics_from_450}, \ref{SAMSARA_numerics_from_450}
and \ref{SAMSARAvsLDFP_convergence_from_450}.

\section{Conclusion and Future Prospects}

Although this is a time dependent inverse problem, we have
considered that we receive certain number of measurements
at a fixed time instant. However, we still aim to observe
expected degradation in the convergence rates in the pre-image
space based on different number of the measurements and the noise amount.
In conclusion, more observations in the image space imply 
better convergence rate in the pre-image space.
As expected from any inverse ill-posed problems, it also has been
observed that less amount of noise in the image space
implies better convergence rate in the pre-image space.

Due to the physical property of the targeted medium,
the actual task is reconstruction of some non-negative
function $\varphi : \Omega \subset \R^{3} \rightarrow \R_{+}.$
However, we have formulated an unconstrained, smooth, convex minimization
problem (\ref{problem}). When the non-negativity constraint is taken
into account, a new minimization problem must be stated as such,

\bea
\label{const_problem}
\varphi_{\alpha}^{\delta} \in \argmin_{\varphi \in \cV} 
\left\{ \frac{1}{2} \norm{\cT \varphi - f^{\delta}}_{\cH}^2 + \alpha J(\varphi) \right\} ,
\\
\varphi(x) \geq 0, \mbox{ for } x \in \Omega \subset \Omega_{o} ,
\eea
where $J : \cV \rightarrow \R_{+}$ is some appropriately
chosen penalty term.

\section{Acknowledgement}

This work was supported by German Weather Service 
(Deutscher Wetterdienst - DWD) project " GPS-Tomography 
for atmospheric data assimilation'' and SFB-755 by University of
G\"{o}ttingen. Author is thankful to D. Russell Luke for the access
to the software SAMSARA for the correct implementation and
the design of the three dimensional total variation penalty term.
Furthermore, author also wishes to thank Thorsten Hohage
for the fruitful discussions over the modelling of the tomography
problem.

%
%
%

\newpage
\bigskip
\section*{References}


 \bibliographystyle{alpha}

\end{document}